\documentclass[12pt]{amsart}
\usepackage{amsmath,amscd,amssymb,amsfonts,graphics}
\usepackage{hyperref}
\newcommand{\hl}{\hyperlink}
\newcommand{\htt}{\hypertarget}
\newcommand{\h}{\hbox}
\newcommand{\q}{\quad}

\newcommand{\bs}{\par\bigskip}
\newcommand{\ms}{\par\medskip}
\newcommand{\sk}{\par\smallskip}
\newcommand{\bsn}{\par\bigskip\noindent}
\newcommand{\msn}{\par\medskip\noindent}

\newcommand{\ges}{\geqslant}
\newcommand{\les}{\leqslant}
\newcommand{\1}{\hskip1pt}
\newcommand{\mopl}{\hbox{$\bigoplus$}}

\newcommand{\msum}{\hbox{$\sum$}}

\newcommand{\D}{{\mathcal D}}
\newcommand{\F}{{\mathcal F}}
\newcommand{\G}{{\mathcal G}}
\newcommand{\Hc}{{\mathcal H}}
\newcommand{\K}{{\mathcal K}}
\newcommand{\Lc}{{\mathcal L}}
\newcommand{\M}{{\mathcal M}}
\newcommand{\OO}{{\mathcal O}}

\newcommand{\Xt}{{}\,\widetilde{\!X}{}}

\newcommand{\sit}{\widetilde{\si}}

\newcommand{\DD}{{\mathbb D}}
\newcommand{\FF}{{\mathbb F}}
\newcommand{\PP}{{\mathbb P}}
\newcommand{\R}{{\mathbb R}}
\newcommand{\RR}{{\mathbf R}}

\newcommand{\Q}{{\mathbb Q}}
\newcommand{\C}{{\mathbb C}}
\newcommand{\N}{{\mathbb N}}
\newcommand{\Z}{{\mathbb Z}}
\newcommand{\V}{{\mathbb V}}

\newcommand{\al}{\alpha}
\newcommand{\ep}{\varepsilon}

\newcommand{\Gr}{{\rm Gr}}
\newcommand{\Om}{\Omega}
\newcommand{\si}{\sigma}
\newcommand{\Si}{\Sigma}
\newcommand{\pp}{{}^{\bf p}}
\newcommand{\MH}{{\rm MH}}
\newcommand{\MHM}{{\rm MHM}}
\newcommand{\Hdg}{{\rm Hdg}}
\newcommand{\Pol}{{\rm Pol}}
\newcommand{\bl}{\bigl}
\newcommand{\br}{\bigr}
\newcommand{\ssb}{\raise.15ex\h{${\scriptscriptstyle\bullet}$}}
\newcommand{\ssc}{\,\raise.15ex\hbox{${\scriptstyle\circ}$}\,}
\newcommand{\onto}{\twoheadrightarrow}
\newcommand{\into}{\hookrightarrow}
\newcommand{\simto}{\,\,\rlap{\hskip1.3mm\raise1.4mm\hbox{$\sim$}}\hbox{$\longrightarrow$}\,\,}
\newcommand{\pl}{\1{+}\1 }
\newcommand{\mi}{\1{-}\1}
\newcommand{\eq}{\,{=}\,}
\newcommand{\gess}{\,{\ges}\,}

\newcommand{\sst}{\,{\subset}\,}

\begin{document}
\title{Hodge modules and cobordism classes}
\author{J. Fern\'andez de Bobadilla}
\address{Javier Fern\'andez de Bobadilla:  
(1) IKERBASQUE, Basque Foundation for Science, Maria Diaz de Haro 3, 48013, 
    Bilbao, Basque Country, Spain;
(2) BCAM,  Basque Center for Applied Mathematics, Mazarredo 14, 48009 Bilbao, 
Basque Country, Spain; 
(3) Academic Colaborator at UPV/EHU.} 
\email{jbobadilla@bcamath.org}
\author{I. Pallar\'es}
\address{Irma Pallar\'es Torres: KU Leuven, Department of Mathematics, Celestijnenlaan 200B, 3001 Leuven, Belgium} 
\email{irma.pallarestorres@kuleuven.be}
\author{M. Saito}
\address{Morihiko Saito: RIMS Kyoto University, Kyoto 606-8502 Japan}
\email{msaito@kurims.kyoto-u.ac.jp}
\thanks{J.F.B. was supported by ERCEA 615655 NMST Consolidator Grant, MINECO by the project 
reference MTM2016-76868-C2-1-P (UCM), by the Basque Government through the BERC 2018-2021 program and Gobierno Vasco Grant IT1094-16, by the Spanish Ministry of Science, Innovation and Universities: BCAM Severo Ochoa accreditation SEV-2017-0718, and by VIASM, through a funded research visit. 
I.P. was supported by ERCEA 615655 NMST Consolidator Grant, by the Basque Government through the BERC 2018-2021 program and Gobierno Vasco Grant IT1094-16, by the Spanish Ministry of Science, Innovation and Universities: BCAM Severo Ochoa accreditation SEV-2017-0718.
M.S was supported by JSPS Kakenhi 15K04816.}
\begin{abstract} We show that the cobordism class of a polarization of Hodge module defines a natural transformation from the Grothendieck group of Hodge modules to the cobordism group of self-dual bounded complexes with real coefficients and constructible cohomology sheaves in a compatible way with pushforward by proper morphisms. This implies a new proof of the well-definedness of the natural transformation from the Grothendieck group of varieties over a given variety to the above cobordism group (with real coefficients). As a corollary, we get a slight extension of a conjecture of Brasselet, Sch\"urmann and Yokura, showing that in the $\Q$-homologically isolated singularity case, the homology $L$-class which is the specialization of the Hirzebruch class coincides with the intersection complex $L$-class defined by Goresky, MacPherson, and others if and only if the sum of the reduced modified Euler-Hodge signatures of the stalks of the shifted intersection complex vanishes. Here Hodge signature uses a polarization of Hodge structure, and it does not seem easy to define it by a purely topological method.
\end{abstract}
\maketitle
\centerline{\bf Introduction}
\msn
F.~Hirzebruch \cite{Hi} introduced the $\chi_y$-characteristic of a compact complex manifold $X$, which is defined as
$$\chi_y(X):=\msum_{p,q\in\N}\,(-1)^q\dim H^q(X,\Om_X^p)\1y^p\,\in\,\Z[y].$$
For $y=-1,0,1$, this specializes respectively to the Euler characteristic, the arithmetic genus, and in the even-dimensional case, the signature of the middle cohomology of $X$. It is the highest degree part of the {\it cohomology Hirzebruch characteristic class\1} $T_{y*}(TX)\in H^{\ssb}(X,\Q)[y]$ of the tangent bundle $TX$. The latter class specializes to the total Chern class $c^*(TX)$, the total Todd class $td^*(TX)$, and the total Thom-Hirzebruch $L$-class $L^*(TX)$ respectively for $y=-1,0,1$.
\sk
This theory has been generalized to the case of singular complex algebraic varieties in \cite{BSY}. We have a natural transformation from the Grothendieck group of varieties over a complex algebraic variety $X$ to the {\it even degree\1} Borel-Moore homology tensored with $\Q[y]$\,:
$$T_{y*}:K_0({\rm Var}/X)\to H^{\rm BM}_{2\1\ssb}(X,\Q)[y].$$
This is compatible with the pushforward by proper morphisms.
\sk
For $y\eq{-}1$, $\,T_{-1*}$ gives the scalar extension of MacPherson's Chern class transformation $c_*{\otimes}\Q$ (see \cite{Ma}) from the constructible functions on $X$. It is shown in \cite{BS} that $c_*(X)$ can be identified with the Schwartz class \cite{Schw} using Alexander duality, when $X$ is compact and embeddable into a complex manifold.
\sk
For $y\eq 0$, $T_{0*}$ is closely related by construction to the Todd class transformation $td_*$ in \cite{BFM}, although $T_{0*}([X])$ does not necessarily coincides with $td_*(X):=td_*([\OO_X])$ unless $X$ has only du Bois singularities.
\sk
For $y\eq 1$, we have the following commutative diagram of abelian groups, assuming $X$ is compact:
\htt{1}{}
$$\begin{array}{rcl}K_0({\rm Var}/X)&\!\!\!\buildrel{sd\,}\over\longrightarrow&\!\!\Om(X)\\ &\!\!\!\!\!\!\!\!\!\!\!\!\raise-1mm\h{$\scriptstyle{T_{1*}}$}\!\!\!\searrow&\,\,\downarrow\,\scriptstyle{L_*}\raise5mm\h{}\\&&\!\!\!\!\!\!\!\!H_{2\1\ssb}(X,\Q)\raise5mm\h{}\end{array}
\leqno(1)$$
Here $\Om(X)$ is the cobordism class group of self-dual bounded $\Q$-complexes with constructible cohomology sheaves on $X$ (which are direct sums of symmetric and skew-symmetric ones), and $L_*$ is the homology $L$-class transformation, see \cite{CS}, \cite{Ba1} (and also \cite{Yo}, \cite{BSY}, where the problem of ambiguities of mapping cones is treated).
The horizontal morphism $sd\1$ is defined in \cite{BSY} using a non-trivial theorem of F.~Bittner \cite{Bi} on the structure of $K_0({\rm Var}/X)$.
The commutativity of (\hl{1}{1}) can be reduced to the assertion on the specialization at $y\eq 1$ in the $X$ smooth compact case explained above, since we have the compatibility of the morphisms in (\hl{1}{1}) with the pushforward by projective morphisms. (Note that $K_0({\rm Var}/X)$ is generated by classes of smooth varieties which are {\it projective\1} over $X$, using affine stratifications and smooth relative projective compactifications of strata such that the divisors at infinity are divisors with simple normal crossings.)
\sk
As is remarked in \cite{BSY}, the image of $[X]$ by the morphism $sd$ does not necessarily coincide with the cobordism class of the intersection complex $[{\rm IC}_X\Q]_{\Om}$, and $T_{1*}([X])$ can be different from the intersection complex $L$-class $L_*([{\rm IC}_X\Q]_{\Om})$ for an irreducible variety $X$ in general.
Note that $L_*([{\rm IC}_X\Q]_{\Om})$ coincides with the homology $L$-class constructed in \cite{GM} as is noted in \cite{BSY}.
It is conjectured in the end of \cite[Remark 0.1]{BSY} that we have a coincidence for the homology classes if $X$ is a compact $\Q$-{\it homology manifold.}
\sk
We extend this to the case where the singularities of $X$ are $\Q$-{\it homologically isolated,} that is, $X\,{\setminus}\,\Si$ is a $\Q$-homology manifold with $\Si\subset X$ a finite subset. The {\it reduced modified Euler-Hodge signature\1} of the stalk of the shifted intersection complex $\mopl_i({\rm IC}_{X_i}\Q[-d_{X_i}])_x$ is defined by
$$\sit_x:=\msum_{j\in\Z}(-1)^j\,{\si}_x^j-1\q\q(x\in X).$$
Here the $X_i$ are irreducible components of $X$ (with $d_{X_i}:=\dim X_i$), and $\si_x^j$ is the sum of the modified Hodge signatures $\si_x^{j,2k}$ of $\Gr^W_{2k}H^j_x$ over $k\,{\in}\,\Z$ with
$$H^j_x:=\mopl_{X_i\ni x}\,\Hc^j({\rm IC}_{X_i}\Q[-d_{X_i}])_x\,.$$
We mean by {\it modified Hodge signature\1} the signature of a polarization of Hodge structure {\it modified as in $(\hl{4}{4})$ below,} that is,
$$\aligned&\si_x^{j,2k}:=\msum_{p\in\Z}\,(-1)^p\1h^{j,\,p,2k-p}_x\\&\q\h{with}\q h^{j,\,p,2k-p}_x:=\dim\Gr^p_F\Gr^W_{2k}(H^j_x)_{\C}.\endaligned$$
Here $(-1)^{k-p}\eq i\1^{q-p}\eq i\1^{p-q}$ if $p\pl q\eq 2k$, and $(-1)^k\eq(-1)^{w(w+1)/2}$ if $w\eq 2k$, see also a remark after Theorem~\hl{T3}{3} below. Note that $\,\si_x^{j,k}\eq 0$ if $k$ is odd, and $\sit_x\eq 0$ if $x\in X\setminus\Si$. (It does not seem easy to define Hodge signatures by a purely topological method.)
\sk
We have the following (which is also conjectured in the end of \cite[Remark 0.1]{BSY}).
\par\htt{T1}{}\msn
{\bf Theorem~1.} {\it If $X$ is a complex variety having only $\Q$-homologically isolated singularities and $\,\sit_x=0$ for any $x\in\Si$, then we have the equality}
\htt{2}{}
$$sd([X])_{\R}=\msum_i\,[{\rm IC}_{X_i}\R]_{\Om}\q\h{in}\,\,\,\,\,\Om\1_{\R}(X).
\leqno(2)$$
\ms
Here $\Om\1_{\R}(X)$ denotes the cobordism class group of self-dual bounded $\R$-complexes with constructible cohomology sheaves on $X$ (which are direct sums of symmetric and skew-symmetric ones), see \hl{1.1}{1.1} below. The left-hand side of (\hl{2}{2}) is defined by the scalar extension of self-dual $\Q$-complexes by $\Q\,{\into}\,\R$. (Note that $\Om\1_{\R}(X)\ne\Om(X)\otimes\R$.) Theorem~\hl{T1}{1} is {\it false\1} in the $\Q$-coefficient case (even if $X$ is a $\Q$-homology manifold), see Proposition~\hl{P2}{2} below.
\ms
Using the commutative diagram (\hl{1}{1}) and the argument as in the proof of Theorem~\hl{T1}{1}, we can show the following.
\par\htt{T2}{}\msn
{\bf Theorem~2.} {\it If $X$ is a compact connected variety having only $\Q$-homologically isolated singularities, then}
\htt{3}{}
$$\msum_{x\in\Si}\,\sit_x\eq 0\,\iff\,T_{1*}(X)\eq\msum_i\,L_*([{\rm IC}_{X_i}\Q]_{\Om})\,\,\,\h{\it in}\,\,H_{2\1\ssb}(X,\Q).
\leqno(3)$$
\ms
The hypothesis of Theorem~\hl{T1}{1} and the left-hand side of (\hl{3}{3}) in Theorem~\hl{T2}{2} hold if the $d_{X_i}$ are {\it even\1} and $X$ has only {\it isolated hypersurface\1} singularities with {\it semisimple\1} Milnor monodromies (for instance, weighted homogeneous isolated hypersurface singularities), see Remark~\hl{R2.3a}{2.3a} below. The above conditions are satisfied in certain other cases, for instance, if $X_1\,{\cap}\,X_2\eq\{x_0\}$ and $X_1,X_2$ have unique singular points at $x_0$ which are ordinary double points with $d_{X_1}$, $d_{X_2}$, $(d_{X_1}{-}d_{X_2})/2$ odd. There are examples such that the left-hand side of (\hl{3}{3}) holds with $X$ {\it globally irreducible,} see Remark~\hl{R2.3b}{2.3b}-\hl{R2.3d}{d} below. There is also an example with $\sit_x\eq0$, $d_X$ {\it odd,} and $(X,x)$ {\it analytic-locally} \h{\it irreducible\1}; for instance, in the case $X$ is defined by $f\eq x^4{+}y^5{+}z^6{+}w^7{+}xyzw$, where the unipotent monodromy part of the vanishing cohomology has dimension 4 with two Jordan blocks of size 3 and 1 (since GCD$(4,6)\ne 1$). This follows from \cite{JKSY} together with \cite{des}. One can also employ a computer program like Singular \cite{DGPS}, see Remark~\hl{R2.3b}{2.3b} below. (It is also possible to apply the Thom-Sebastiani theorem to $T_{p_1,p_2,p_3}$ and $T_{q_1,q_2,q_3}$ with $p_i,q_j$ mutually prime for any $i,j$.)
\sk
In the $\Q$-homology manifold case (that is, if $\Si=\emptyset$), Theorem~\hl{T2}{2} implies the conjecture in \cite{BSY} mentioned above. This has been proved in special cases by many people; for instance, in the isolated hypersurface singularity case \cite{CMSS1}, quotient singularity case \cite{CMSS2}, some toric variety case \cite{MS}, some threefold case \cite{Ba2}. The conjecture is shown in \cite{FP} in the $X$ projective case (although the argument there, especially the proof of Proposition~\hl{P1}{1} below, does not seem necessarily easy to follow for everybody except experts in representation theory). It is possible to remove the projectivity assumption in \cite{FP} by using \cite{mhc}.
\sk
For the proofs of Theorems~\hl{T1}{1} and \hl{T2}{2}, we need the following
\par\htt{P1}{}\msn
{\bf Proposition~1.} {\it Let $\M$ be a pure $A$-Hodge module of weight $w$ on $X$ with $A$ a subfield of $\R$. Let $K_{\R}^{\ssb}$ be its underlying $\R$-complex $($that is, $K_{\R}^{\ssb}\eq K_A^{\ssb}{\otimes}_A\R$ with $K_A^{\ssb}$ its underlying $A$-complex$)$. Let
$$S_{\R}:K_{\R}^{\ssb}{\otimes}K_{\R}^{\ssb}\to\DD\R_X(-w)$$
be the scalar extension of a polarization $S$ of the $A$-Hodge module $\M$, where $\DD\R_X$ denotes the dualizing complex. Then the cobordism class $[K_{\R}^{\ssb},S_{\R}]\in\Om\1_{\R}(X)$ does not depend on the choice of a polarization $S$.}
\ms
A stronger assertion is used in \cite{FP}, where the assumption that $K^{\ssb}_{\R}$ underlies a pure $\R$-Hodge module and the pairing gives a polarization of it is replaced by the conditions that $K^{\ssb}_{\R}$ is an intersection complex with coefficients in a local system, the stalk of the local system at some point underlies a Hodge structure, and the restriction of the pairing gives a polarization of it. This stronger claim can be shown by the same argument as in \hl{2.1}{2.1} below. The argument in that paper employs a highly sophisticated method from representation theory, and demonstrates a much stronger assertion (see Remark~\hl{R2.1}{2.1} below), where the {\it semi-simplicity\1} of local system is {\it not\1} needed as long as we have a smooth 1-parameter family of non-degenerate pairings. (This is quite surprising.) In this paper we give a quite simple proof of Proposition~\hl{P1}{1} using Hodge theory, see \hl{2.1}{2.1} below.
\sk
From Proposition~\hl{P1}{1}, we can deduce the following.
\par\htt{T3}{}\msn
{\bf Theorem~3.} {\it For a subfield $A\subset\R$, there is a natural transformation
\htt{X}{}
$$\Pol:K_0\bl(\MHM(X,A)\br)\to\Om\1_{\R}(X)$$
defined by
\htt{4}{}
$$\Pol([\M])=[K_{\R}^{\ssb},(-1)^{w(w+1)/2}S_{\R}]\,\,\,\in\,\,\,\Om\1_{\R}(X),
\leqno(4)$$
for a pure $A$-Hodge module $\M$ of weight $w$ on $X$, where $K_{\R}^{\ssb}$ and $S_{\R}$ are as in Proposition~$\hl{P1}{1}$. This is compatible with the pushforward by a proper morphisms $f:X\to Y$, that is,}
\htt{5}{}
$$\Pol\ssc f_*=f_*\ssc\Pol.
\leqno(5)$$
\ms
The sign $(-1)^{w(w+1)/2}$ in (\hl{4}{4}) is closely related to the sign $(-1)^{d_Z(d_Z-1)/2}$ in Theorem~\hl{T1.3}{1.3} below. We get an additional sign $(-1)^{d_Z}$ from the shift of complex by $d_Z$ in the constant coefficient case. This explains the reason for which we need no sign in (\hl{2}{2}). The sign in (\hl{4}{4}) with $X\eq{\rm pt}$ is also compatible with a formula for $\chi_1$ of polarized Hodge structures of even weights written after (2-14) in \cite{Schu} (as is remarked by the referee). Here $\chi_1(V)$ is the specialization at $y\eq 1$ of $\chi_y(V):=\msum_{p,q}\,h^{p,q}(V)(-y)^p$ for a mixed Hodge structure $V$. For a compact pure-dimensional variety $Z$, the last assertion of Theorem~\hl{T3}{3} then implies the {\it Hodge index theorem\1} for intersection cohomology (as is mentioned by the referee):
\htt{6}{}
$$\chi_1\bl({\rm IH}^{\ssb}(Z)\br)={\rm sign}\bl({\rm IH}^{d_Z}(Z)\br).
\leqno(6)$$
(This follows also from \cite[Thm.\,5.3.1]{mhp}, that is, Theorem~\hl{T1.4}{1.4} below.)
\sk
From Theorem~\hl{T3}{3} we can deduce the following.
\par\htt{C1}{}\msn
{\bf Corollary~1.} {\it We have the commutative diagram}
\htt{7}{}
$$\begin{array}{rcl}K_0({\rm Var}/X)&\buildrel{\!\!\!\Hdg}\over\longrightarrow&\!\!K_0\bl(\MHM(X,A)\br)\\ &\!\!\!\!\!\!\!\!\raise-1mm\h{$\scriptstyle{sd\hskip.5pt{}_{\R}}$}\!\!\!\searrow&\,\,\,\,\,\downarrow\,\scriptstyle{\Pol}\raise5mm\h{}\\&&\!\Om\1_{\R}(X)\raise5mm\h{}\end{array}
\leqno(7)$$
\ms
Recall that the horizontal morphism $\Hdg$ is defined by
$$[f\,{:}\,Y\1{\to}\1X]\,\mapsto\,\msum_{j\in\Z}\,(-1)^j[H^jf_!\1A_{h,Y}],$$
see \cite{BSY}. (Its well-definedness follows from (\hl{1.4.5}{1.4.5}) below.) Here $A_{h,Y}\in D^b\bl(\MHM(Y,A)\br)$ is defined by $a_Y^*A_h$ with $a_Y:Y\to{\rm pt}$ the structure morphism and $A_h$ the $A$-Hodge structure of rank $1$ with weight $0$, see \hl{1.4}{1.4} below. In the case $f\eq{\rm id}$ and $X$ is smooth, we have the isomorphism
$$A_{h,X}=(A_{h,X}[d_X])[-d_X],$$
where $A_{h,X}[d_X]$ is a pure $A$-Hodge module of weight $d_X\,(:=\dim X)$. So the commutativity of (\hl{7}{7}) holds for $[X]\in K_0({\rm Var}/X)$, since a canonical polarization of $A_{h,X}[d_X]$ is defined with the sign $(-1)^{d_X(d_X-1)/2}$ (see \cite[5.2.12 and 5.4.1]{mhp}) and $d_X(d_X{-}1)/2\pl d_X\eq d_X(d_X{+}1)/2$. The general case then follows from the compatibility with the pushforward by projective morphisms.
\sk
The proof of the commutativity of $\Pol$ with the pushforward by {\it projective\1} morphisms is similar to the argument for the coincidence of $\chi_1(X)$ with the signature of the middle cohomology in the $X$ smooth projective case, see \cite[Thm.\,15.8.2]{Hi} (and also \cite[3.6]{MSS}). Here we use the property of sign noted in (\hl{2.2.5}{2.2.5}) below, and we need in an essential way the hard Lefschetz property and the positivity of the induced polarization on the primitive part up to a given sign (see \cite[Thm.~5.3.1]{mhp} and also Theorem~\hl{T1.4}{1.4} below) together with a certain cobordism relation (see \cite{CS}, \cite{Yo} and also Proposition~\hl{P1.1a}{1.1a} below). The proper morphism case can be reduced to the projective morphism case by induction on dimension.
\sk
The above argument implies the well-definedness of $sd\1_{\R}$ (only in the $\R$-coefficient case), which was proved in \cite{BSY} using a highly non-trivial theorem of F.~Bittner \cite{Bi}.
Theorem~\hl{T1}{1} then follows from the commutativity of (\hl{7}{7}) in Corollary~\hl{C1}{1} using the canonical morphism
\htt{8}{}
$$\Q_{h,X}\to\mopl_i\,{\rm IC}_{h,X_i}\Q[-d_{X_i}]\q\h{in}\,\,\,\,D^b\MHM(X,\Q).
\leqno(8)$$
This morphism is obtained by the adjunction property of $a_X^*$ and $(a_X)_*$ with $a_X:X\to{\rm pt}$ the structure morphism. (It is an isomorphism if and only if $X$ is a $\Q$-homology manifold, see for instance \cite[p.~34]{BM}.) For details, see \hl{2.3}{2.3} below.
\sk
Considering the mapping cone of (\hl{8}{8}), we can also get the following equalities for a compact variety $X$ (as is pointed out by the referee):
\htt{9}{}
$$T_{1,k}(X)=\msum_i\,L_k([{\rm IC}_{X_i}\Q]_{\Om})\,\,\,\,\h{in}\,\,\,\,H_{2k}(X,\Q)\q\h{if}\q k>\dim\Si.
\leqno(9)$$
Here $\Si\into X$ is the smallest closed subset such that $X\,{\setminus}\,\Si$ is a $\Q$-homology manifold.
\sk
Finally note that we have the following.
\par\htt{P2}{}\msn
{\bf Proposition~2.} {\it Theorem~$\hl{T1}{1}$ does not hold with $\Q$-coefficients.}
\ms
This can be shown by using classical Witt groups \cite{MH}, see \hl{2.5}{2.5} below.
\sk
In Section~\hl{S1}{1} we review some basics of cobordism classes and Hodge modules.
In Section~\hl{S2}{2} we prove Theorems~\hl{T1}{1}--\hl{T3}{3} after showing Proposition~\hl{P1}{1}.
\sk
We thank J.~Sch\"urmann for a good question which inspired us with Proposition~\hl{P1.1d}{1.1d} and Theorem~\hl{T1.1}{1.1} below. We are grateful to the referee for very careful reading of the paper.
\par\htt{i}{}\msn
{\bf Conventions.}~(i) In this paper, a variety means effectively the associated analytic space of a reduced separated scheme of finite type over $\C$. (We need the classical topology, since local systems and more generally constructible sheaves are used, although algebraic coherent sheaves are not.) However, the algebraic structure is not forgotten, since morphisms of varieties are always induced by morphisms of schemes; in particular, a subvariety means the associated analytic space of a reduced subscheme, and stratifications are always algebraic.
\par\htt{ii}{}\ms
(ii) In the definition of polarization of Hodge structure, we have the action of the Weil operator on the second factor as in \cite{De1}. This produces the difference of sign by $(-1)^w$ from the other convention, where $w$ is the weight of Hodge structure. This and the sign $(-1)^{d_Z(d_Z-1)/2}$ in Theorem~\hl{T1.3}{1.3} below give a ``dictionary" between the sign system in this paper and that in \cite{FP}.
\par\htt{iii}{}\ms
(iii) We denote by $\pp\tau_{\les j}$ and $\pp\Hc^j$ the truncation and cohomology functor for the middle perversity, which is constructed in \cite{BBD}. We denote by $D^b_c(X,A)^{[k]}$ the full subcategory of $ D^b_c(X,A)$ defined by the condition that $\pp\Hc^jK\eq 0$ ($j\ne k$) for $K\in D^b_c(X,A)$, where $k\in\Z$ and $A$ is a subfield of $\C$. These are abelian categories, and $D^b_c(X,A)^{[0]}$ is stable by the functor $\DD$ associating the dual, see \cite{BBD}.
\bs\htt{S1}{}\bs
\vbox{\centerline{\bf 1. Preliminaries}
\bsn
In this section we review some basics of cobordism classes and Hodge modules.}
\par\htt{1.1}{}\msn
{\bf 1.1.~Cobordism classes.} Let $X$ be a complex algebraic variety. Let $A$ be any field of characteristic 0 in this subsection. Let $(\F,S)$ be a self-dual bounded $A$-complexes with constructible cohomology sheaves. This is called a {\it self-dual complex\1} for short. More precisely, it means that $\F\in D^b_c(X,A)$ and $\F$ is endowed with a perfect pairing
$$S:\F\otimes\F\to\DD A_X,$$
where $\DD A_X$ is the dualizing complex with $\DD$ the functor associating the dual. (Here Tate twists are omitted to simplify the notation, choosing $\sqrt{-1}\in\C$ if necessary, although there is no canonical choice.) In general a morphism
$$\F\otimes\G\to\DD A_X$$
is called a {\it perfect pairing\1} if the corresponding morphism $\F\to\DD\G$ is an isomorphism in $D^b_c(X,A)$.
Here we use the following canonical isomorphism for $\F,\G,\K\in C^b_c(X,A)$\,:
\htt{1.1.1}{}
$${\rm Hom}(\F\,{\otimes}\,\G,\K)={\rm Hom}\bl(\F,\Hc om(\G,\K)\br).
\leqno(1.1.1)$$
We may assume that $\K$ is an injective complex taking a resolution (where $\K\in C^+_c(X,A)$). Then (\hl{1.1.1}{1.1.1}) holds in $D^b_c(X,A)$, see \cite{BBD}, etc. Recall that $\DD$ is defined by
$$\DD\F:=\Hc om(\F,a_X^!A_{\rm pt}),$$
with $a_X:X\to{\rm pt}$ the natural morphism, hence $\DD A_X=a_X^!A_{\rm pt}$.
\sk
We say that the pairing $S$ is {\it symmetric\1} (or \h{\it skew-symmetric\1}) if $S\ssc\iota\eq S$ (or $-S$). Here we denote by $\iota:\F{\otimes}\G\to\G{\otimes}\F$ the involution defined by
\htt{1.1.2}{}
$$\F^i{\otimes}\G^j\,\ni\,u{\otimes}v\,\mapsto\,(-1)^{ij}v{\otimes}u\,\in\,\G^j{\otimes}\F^i,
\leqno(1.1.2)$$
see \cite{De2}, \cite{De3}, \cite[1.3.4]{MSS}, etc.
With this definition it is easy to see that the canonical self-pairing on $A_X[d_X]$ is $(-1)^{d_X}$-symmetric (that is, symmetric if $d_X$ is even, and skew-symmetric otherwise) when $X$ is smooth and pure-dimensional, since $\DD A_X\cong A_X[2d_X]$ (here the isomorphism may depend on the choice of $\sqrt{-1}$). The argument seems more complicated if one adopts a definition of (skew-)symmetric pairing using the functor $\DD$, see Remark~\hl{R1.1e}{1.1e} below.
\sk
We assume $S$ is either symmetric or skew-symmetric. The {\it symmetric cobordism group\1} $\Om_{A+}(X)$ is defined to be the quotient of the monoid of isomorphism classes of symmetric self-dual bounded $A$-complexes with constructible cohomology sheaves, which is divided by the {\it cobordism relation\1} explained below. The {\it skew-symmetric cobordism group\1} $\Om_{A-}(X)$ is defined by replacing symmetric with skew-symmetric. Here the sum is given by direct sum. These are abelian groups, since we have in $\Om_{A+}(X)$ or $\Om_{A-}(X)$
\htt{1.1.3}{}
$$[(\F,S)]+[(\F,-S)]=0,
\leqno(1.1.3)$$
see Proposition~\hl{P1.1b}{1.1b} below (where $\F\eq\F'\eq\G\eq\G'$). Finally set
$$\Om_A(X):=\Om_{A+}(X)\oplus\Om_{A-}(X).$$
Sometimes $A$ is omitted when $A=\Q$.
\sk
We say that $(\F,S)$ is {\it cobordant to\1} $(\F',S')$ if there are $(\F_i,S_i)$ $(i\in[0,r])$ such that $(\F_0,S_0)=(\F,S)$, $(\F_r,S_r)=(\F',S')$, and $(\F_{i-1},S_{i-1})$ is directly cobordant to $(\F_i,S_i)$ for any $i\in[1,r]$ (see the definition just below), where the $(\F_i,S_i)$ are either symmetric for any $i$ or skew-symmetric for any $i$.
\sk
We say that $(\F,S)$ {\it directly cobordant to\1} $(\F',S')$ if there is a commutative diagram
\htt{1.1.4}{}
$$\begin{array}{ccc}\G&\buildrel{\rho'}\over\to&\F'\\ \!\!\!{\scriptstyle\pi}\!\downarrow&\raise4mm\h{}\raise-2mm\h{}&\,\,\downarrow\!{\scriptstyle\pi'}\\\F&\buildrel{\rho}\over\to&\G'\end{array}
\leqno(1.1.4)$$
in $D^b_c(X,A)$ together with a perfect pairing $S'':\G\otimes\G'\to\DD A_X$ such that
\htt{1.1.5}{}
$$\aligned S\ssc(\pi\1{\otimes}\1{\rm id})&=S''\ssc({\rm id}\1{\otimes}\1\rho):\G\otimes\F\to\DD A_X,\\ S'\ssc(\rho'\1{\otimes}\1{\rm id})&=S''\ssc({\rm id}\1{\otimes}\1\pi'):\G\otimes\F'\to\DD A_X,\endaligned
\leqno(1.1.5)$$
and moreover the morphism of mapping cones $C(\rho')\to\C(\rho)$ induced (non-canonically) by $(\pi,\pi')$ is an isomorphism in $D^b_c(X,A)$.
\sk
In the case $\F,\F'$ belong to $D^b_c(X,A)^{[0]}$ (see Convention (\hl{iii}{iii}) at the end of the introduction), we say that $(\F,S)$ is {\it directly subquotient cobordant to\1} $(\F',S')$ if the above conditions are satisfied with $\G,\G'\in D^b_c(X,A)^{[0]}$ and moreover the horizontal and vertical morphisms of (\hl{1.1.4}{1.1.4}) are either {\it injective\1} and {\it surjective} respectively or {\it surjective\1} and {\it injective} respectively. We say that $(\F,S)$ is {\it subquotient cobordant\1} to $(\F',S)$ if there are $(\F_i,S_i)$ $(i\in[0,r])$ such that $(\F_0,S_0)=(\F,S)$, $(\F_r,S_r)=(\F',S')$, and $(\F_{i-1},S_{i-1})$ is directly subquotient cobordant to $(\F_i,S_i)$ for any $i\in[1,r]$. Here we assume $\F,\F',\F_i\in D^b_c(X,A)^{[0]}$
\sk
We will see that cobordism relation is essentially equivalent to subquotient cobordism relation, see Theorem~\hl{T1.1}{1.1} below.
\par\htt{R1.1a}{}\msn
{\bf Remark~1.1a.} Since the mapping cone is {\it unique up to a non-canonical isomorphism,} the isomorphism condition for mapping cones after (\hl{1.1.5}{1.1.5}) is well-defined. Moreover this condition is equivalent to a similar condition for $C(\pi)\to C(\pi')$. Indeed, the diagram (\hl{1.1.4}{1.1.4}) is commutative up to a homotopy $h$ at the level of complex, and we can construct a double complex representing the mapping cone of the mapping cone using this homotopy $h$ as is well-known in the theory of derived categories. The condition after (\hl{1.1.5}{1.1.5}) is equivalent to the acyclicity of the associated single complex, see also \cite[Prop.~1.1.11]{BBD}.
\sk
The diagram (\hl{1.1.4}{1.1.4}) and the construction of the associated octahedral diagram are closely related to a commutative diagram in the proof of \cite[Prop.~6.7]{Yo}.
There is a problem of the ambiguity of mapping cones in the definition of cobordism class in \cite{CS}. This is improved in \cite{Yo}, see also \cite{BSY}. However, the new definition might impose unnecessary complications. It seems that (\hl{1.1.4}{1.1.4}) is the essential part of the octahedral diagram in the definition of cobordism relation.
\par\htt{R1.1b}{}\msn
{\bf Remark~1.1b.} The condition that $(\F,S)$ is directly cobordant to $(\F',S')$ is {\it symmetric\1} (and also reflexive).
\par\htt{R1.1c}{}\msn
{\bf Remark~1.1c.} Perfect pairings and cobordism relations are stable by the pushforward $f_*$ under a {\it proper\1} morphism of complex varieties $f:X\to Y$. Indeed, if we have a perfect pairing $\F\1{\otimes}\1\G\to\DD A_X$, then we have the induced pairing by
\htt{1.1.6}{}
$$\RR f_*\F\1{\otimes}\1\RR f_*\G\to\RR f_*(\F\1{\otimes}\1\G)\to\RR f_*\DD A_X\buildrel{\!\!\!\rm Tr}\over\to\DD A_Y,
\leqno(1.1.6)$$
and this is also a perfect pairing by Verdier duality. We thus get the pushforward
$$f_*:\Om_A(X)\to\Om_A(Y),$$
compatible with the direct sum decomposition written after (\hl{1.1.3}{1.1.3}).
\ms
The following properties are well-known, see \cite{CS}, \cite{Yo}, \cite{SW}, etc.
\par\htt{P1.1a}{}\msn
{\bf Proposition~1.1a.} {\it A self-dual complex $(\F,S)$ is directly cobordant to $\pp\Hc^0(\F,S)$, hence we have the equality}
\htt{1.1.7}{}
$$[(\F,S)]=[\pp\Hc^0(\F,S)]\q\h{in}\,\,\,\,\Om_A(X).
\leqno(1.1.7)$$
\msn
{\it Proof.} Consider the commutative diagram as in (\hl{1.1.4}{1.1.4}):
$$\begin{array}{ccc}\pp\tau_{\les0}\F&\buildrel{\rho'}\over\to&\F\\ \!\!\!{\scriptstyle\pi}\!\downarrow&\raise4mm\h{}\raise-2mm\h{}&\,\,\downarrow\!{\scriptstyle\pi'}\\ \pp\Hc^0\F&\buildrel{\rho}\over\to&\pp\tau_{\ges0}\F\end{array}$$
see Convention (\hl{iii}{iii}).
Then the perfect pairing $S$ on $\F{\otimes}\F$ induces canonically the perfect pairings on $\pp\Hc^0\F\1{\otimes}\1\pp\Hc^0\F$ and $\pp\tau_{\les0}\F\1{\otimes}\1\pp\tau_{\ges0}\F$ by the theory of $t$-structure in \cite{BBD} via the isomorphism (\hl{1.1.1}{1.1.1}). So Proposition~\hl{P1.1a}{1.1a} follows.
\par\htt{P1.1b}{}\msn
{\bf Proposition~1.1b.} {\it Assume a self-dual complex $(\F,S)$ is directly cobordant to $(\F',S')$. Then the direct sum $(\F',S')\oplus(\F,-S)$ is directly cobordant to $0$, hence}
\htt{1.1.8}{}
$$[(\F',S')]+[(\F,-S)]=0\q\h{\it in}\,\,\,\,\Om_A(X).
\leqno(1.1.8)$$
\msn
{\it Proof.} Let $\G,\G'$ and $\pi,\rho,\pi',\rho'$ be as in (\hl{1.1.4}{1.1.4}) associated to the direct cobordism between $(\F,S)$ and $(\F',S')$. We have the commutative diagram
$$\begin{array}{ccc}\G&\buildrel{(\rho',\pi)}\over\longrightarrow&\!\!(\F',S')\oplus(\F,-S)\\\downarrow&\raise4mm\h{}\raise-2mm\h{}&\q\q\,\downarrow{\scriptstyle(\pi',-\rho')}\\ 0&\q\q\longrightarrow&\!\!\!\G'\end{array}$$
Moreover (\hl{1.1.5}{1.1.5}) and the condition after it are satisfied for this diagram (using Remark~\hl{R1.1a}{1.1a}). So Proposition~\hl{P1.1b}{1.1b} follows.
\par\htt{P1.1c}{}\msn
{\bf Proposition~1.1c.} {\it For a self-dual complex $(\F,S)$, we have the equality
\htt{1.1.9}{}
$$[(\F,S)]=\mopl_{i=1}^r[(\F_i,S_i)]\q\h{\it in}\,\,\,\,\Om_A(X),
\leqno(1.1.9)$$
with $\F_i$ simple objects of $D^b_c(X,A)^{[0]}$, see Convention~{\rm (\hl{iii}{iii})}.}
\ms
(It is known that simple objects of $D^b_c(X,A)^{[0]}$ are intersection complexes with coefficients in simple local systems defined on locally closed irreducible smooth subvarieties, see \cite{BBD}.)
\msn
{\it Proof.} We may assume $\F\in D^b_c(X,A)^{[0]}$ by Proposition~\hl{P1.1b}{1.1b}. We proceed by induction on the length of $\F$ in the abelian category $D^b_c(X,A)^{[0]}$, which is Noetherian and Artinian. If $\F$ is simple, the assertion is trivial. Let $\G$ be a simple subobject of $\F$. Let $\G'$ be its orthogonal complement with respect to the perfect pairing $S$, that is, the kernel of the composition
$$\phi:\F\simto\DD\F\onto\DD\G,$$
where the first isomorphism corresponds to the perfect pairing $S$. The canonical morphism $\G\to\F/\G'$ is identified with the composition of the inclusion $\G\into\F$ with the above composition $\phi$. Since $\G$ is a simple object, it is either an isomorphism or zero.
\sk
If it is an isomorphism, then we have the direct sum decomposition $\F\eq\G\oplus\G'$ compatible with $S$, and we can apply the inductive hypothesis to $\G'$.
\sk
If it is zero, then we have the inclusion $\G\subset\G'$. We see that the condition for direct subquotient cobordism relation in (\hl{1.1.4}{1.1.4}--\hl{1.1.5}{5}) is satisfied for the diagram
$$\begin{array}{ccc}\G'&\into&\F\\ \rlap{\h{$\downarrow$}}\raise1mm\h{$\downarrow$}&\raise4mm\h{}\raise-2mm\h{}&\rlap{\h{$\downarrow$}}\raise1mm\h{$\downarrow$}\\\G'/\G&\into&\F/\G\end{array}$$
so that $\G'/\G$ is directly subquotient cobordant to $\F$. We can then apply the inductive hypothesis to $\G'/\G$. Thus Proposition~\hl{P1.1c}{1.1c} is proved.
\par\htt{P1.1d}{}\msn
{\bf Proposition~1.1d.} {\it If $(\F,S)$ is directly cobordant to $(\F',S')$, then there is $(\F'',S'')$ such that $\F''\in D^b_c(X,A)^{[0]}$ and $(\F'',S'')$ is directly subquotient cobordant to both $\pp\Hc^0(\F,S)$ and $\pp\Hc^0(\F',S')$.}
\msn
{\it Proof.} Setting $\K:=C(\rho)$, $\K':=C(\rho')$, we get a commutative diagram
\htt{1.1.10}{}
$$\begin{array}{ccccc}\pp\Hc^{-1}\K'&\rightarrow&\pp\Hc^0\G&\rightarrow&\pp\Hc^0\F'\\ \raise-.5mm\h{\rotatebox{90}{$\cong$}}&\raise4mm\h{}\raise-1mm\h{}&\,\,\,\downarrow\!\!\raise.7mm\h{$\scriptstyle\pi^0$}&\searrow\!\!\!\!\raise1.5mm\h{$\scriptstyle\delta$}&\downarrow\\ \pp\Hc^{-1}\K&\rightarrow&\pp\Hc^0\F&\buildrel{\rho^0}\over\rightarrow&\pp\Hc^0\G'\end{array}
\leqno(1.1.10)$$
in the abelian category $D^b_c(X,A)^{[0]}$. (This is quite similar to a commutative diagram in the proof of \cite[Prop.~6.7]{Yo}.) Here $\pi^0$ denotes $\pp\Hc^0\pi$, etc., and the horizontal sequences are exact. These imply the inclusion
$$\Lc'':={\rm Ker}\,\rho^0\,\,\subset\,\,\Lc:={\rm Im}\,\pi^0\q\q\h{with}\q\q\Lc/\Lc''\simto{\rm Im}\,\delta,$$
where $\delta:=\rho^0\ssc\pi^0$. Set
$$\Lc':=\pp\Hc^0\F/\Lc''\,(={\rm Coim}\,\rho^0).$$
Since $\rho^0$ is identified with the dual of $\pi^0$ by (\hl{1.1.5}{1.1.5}), we see that $\Lc'$ is identified with the dual of $\Lc\eq{\rm Im}\,\pi^0$ by $S$, and $S$ induces a self-pairing of $\Lc/\Lc''$. Moreover the latter can be identified with the induced self-pairing $S''$ on ${\rm Im}\,\delta\subset\pp\Hc^0\G'$ using (\hl{1.1.5}{1.1.5}). We thus get the commutative diagram as in (\hl{1.1.4}{1.1.4})
$$\begin{array}{ccc}\Lc&\into&\pp\Hc^0\F\\ \rlap{\h{$\downarrow$}}\raise1mm\h{$\downarrow$}&\raise4mm\h{}\raise-2mm\h{}&\rlap{\h{$\downarrow$}}\raise1mm\h{$\downarrow$}\\{\rm Im}\,\delta&\into&\Lc'\end{array}$$
with (\hl{1.1.5}{1.1.5}) satisfied. So $({\rm Im}\,\delta,S'')$ is directly subquotient cobordant to $\pp\Hc^0(\F,S)$. We can apply a similar argument replacing $(\F,S)$ with $(\F',S')$. This finishes the proof of Proposition~\hl{P1.1d}{1.1d}.
\ms
By Propositions~\hl{P1.1a}{1.1a} and \hl{P1.1d}{1.1d}, we get the following.
\par\htt{T1.1}{}\msn
{\bf Theorem~1.1.} {\it A self-dual complex $(\F,S)$ is cobordant to $(\F',S')$ if and only if there are self-dual complexes $(\F_i,S_i)$ with $\F_i\in D^b_c(X,A)^{[0]}$ for $i\in[0,2r]$ such that $(\F_0,S_0)=\pp\Hc^0(\F,S)$, $(\F_{2r},S_{2r})=\pp\Hc^0(\F',S')$, and $(\F_{i-1},S_{i-1})$ is directly subquotient cobordant to $(\F_i,S_i)$ for any $i\in[1,2r]$.}
\par\htt{R1.1d}{}\msn
{\bf Remark~1.1d.} The condition after (\hl{1.1.5}{1.1.5}) implies that the commutative diagram (\hl{1.1.4}{1.1.4}) is completed {\it non-canonically\1} to an octahedral diagram by adding
$$C(\rho')\cong C(\rho)\q\h{and}\q C(\pi)\cong C(\pi')$$
at the vertex in the center of the upper and lower part of the octahedral diagrams respectively, see also \cite{BSY}, \cite{Yo}. Here the octahedral diagram is viewed from a {\it different angle,} since we have a {\it commutative\1} square instead of a {\it circulating\1} one at the boundary. (Note that an octahedral diagram has one circulating boundary square and two commutative boundary squares, where the initial and terminal vertices of each commutative boundary square are the two vertices outside the circulating boundary square. Here a boundary square means that it divides the octahedron into two pyramids.)
We do not demand a self-duality isomorphism of the whole octahedral diagram, since its construction may be rather complicated in general. (Here the problem of sign seems rather non-trivial.)
\sk
Note that there is no problem in the case of {\it subquotient cobordism,} since the \h{\it kernel\1} and \h{\it cokernel\1} in the abelian category $D^b_c(X,A)^{[0]}$ are unique up to unique isomorphism (making certain diagrams commutative).
Here we have an increasing filtration $G$ on $\F'$ in the abelian category $D^b_c(X,A)^{[0]}$ such that $\Gr^G_i\F'\eq0$ ($i\ne-1,0,1$), and
$$\Gr^G_0\F'=\F,\q G_0\F'=\G,\q\F'/G_{-1}\F'=\G'.$$
So the octahedral diagram can be completed {\it canonically\1} by adding $\Gr^G_1\F'$, $\Gr^G_{-1}\F'$ to the commutative diagram
$$\begin{array}{ccc}G_0\F'&\into&\F'\\ \rlap{\h{$\downarrow$}}\raise1mm\h{$\downarrow$}&\raise4mm\h{}\raise-2mm\h{}&\rlap{\h{$\downarrow$}}\raise1mm\h{$\downarrow$}\\ \Gr^G_0\F'&\into&\F'/G_{-1}\F'\end{array}$$
\sk
Note that a short exact sequence in the abelian category $D^b_c(X,A)^{[0]}$ defines an extension class, and the latter should coincide with the connecting morphism in the derived category $D^b_c(X,A)$ (without any sign, see also \cite[2.1]{Ve}).
This would imply that the definition of $\Om(X)$ in this paper is equivalent to the one in \cite{Yo} (admitting that Proposition~\hl{P1.1a}{1.1a} holds with the definition given there, see \cite[Example~6.6]{Yo}).
\par\htt{R1.1e}{}\msn
{\bf Remark~1.1e.} The sign becomes more complicated if we have to deal with a double complex like ${\rm Hom}(\F^{\ssb},\G^{\ssb})$ which is {\it contravariant\1} for the first factor, where we get an additional minus sign, see \cite{De2}. For instance, there is an {\it anti-commutative\1} diagram (see also \cite{BSY})
\htt{1.1.11}{}
$$\begin{array}{ccc}{\rm Hom}(\F^{\ssb}[1],\G^{\ssb}[1])&\simto&{\rm Hom}(\F^{\ssb},\G^{\ssb}[1])[-1]\\ \downarrow\!\!\1\raise.3mm\h{\rotatebox{90}{$\scriptstyle\sim$}}&\raise5mm\h{}\raise-3mm\h{}&\downarrow\!\!\1\raise.3mm\h{\rotatebox{90}{$\scriptstyle\sim$}}\\{\rm Hom}(\F^{\ssb}[1],\G^{\ssb})[1]&\simto&{\rm Hom}(\F^{\ssb},\G^{\ssb})\end{array}
\leqno(1.1.11)$$
This seems to be related closely to the self-duality isomorphism of $A_X[d_X]$, that is,
$$\Hc om(A_X[d_X],A_X[2d_X])=A_X[d_X],$$
where $X$ is smooth. We may have to {\it choose\1} one of the two above isomorphisms, and a consistent sign system may contain a set of consistent choices for all the anti-commutative diagrams. There are also rather complicated isomorphisms like
$$\DD\cong(\DD\ssc\DD)\ssc\DD\cong\DD\ssc(\DD\ssc\DD)\cong\DD,$$
where the determination of signs does not seem quite trivial.
It does not seem very clear whether these problems are completely clarified in the literature. Indeed, there seem to be several sign systems in the literature. (For instance, the signs in the definitions of {\it mapping cone\1} and its {\it associated distinguished triangle\1} do not seem unique. These do not seem to be specified in \cite{CH} where the signs necessary for each isomorphism in (\hl{1.1.11}{1.1.11}) seem to be given, but not the choice of the isomorphism between the top-left and bottom-right of the diagram; there are {\it two choices\1} since we have the {\it anti-commutativity.})
In order to avoid the above problems, it is recommended in \cite{De3} to use perfect pairings rather than the functor $\DD$, and we follow this suggestion in this paper.
\msn
{\bf 1.2.~Relation with Witt groups.} In the notation of \hl{1.1}{1.1}, let $W\!\!_{A+}(X)$ be the monoid of isomorphism classes of symmetric self-dual $A$-complexes divided by the submonoid of isomorphism classes of symmetric self-dual $A$-complexes which are directly cobordant to $0$. More precisely, we have $[(\F,S)]\eq[(\F',S')]$ in $W\!\!_{A+}(X)$ if and only if
$$(\F,S)\pl(\G_1,S_1)\eq(\F',S')\pl(\G_2,S_2),$$
with $(\G_i,S_i)$ directly cobordant to $0$ $(i\eq1,2)$.
(This is compatible with \cite{Bal} by \cite{SW} and Theorem~\hl{T1.2}{1.2} below.)
Define $\pp W\!\!_{A+}(X)$ by assuming further that all the self-dual $A$-complexes belong to $D^b_c(X,A)^{[0]}$. Similarly we can define $W\!\!_{A-}(X)$ and $\pp W\!\!_{A-}(X)$ by replacing symmetric with skew-symmetric. (Note that, if a self-dual complex $(\F,S)$ with $\F\in D^b_c(X,A)^{[0]}$ is directly cobordant to 0, then $(\F,S)$ is directly subquotient cobordant to 0.)
\sk
From Theorem~\hl{T1.1}{1.1} and Propositions~\hl{P1.1a}{1.1a}--\hl{P1.1b}{b}, we can deduce the following (see also \cite{SW}):
\par\htt{T1.2}{}\msn
{\bf Theorem~1.2.} {\it There are canonical isomorphisms of abelian groups}
\htt{1.2.1}{}
$$\aligned\pp W\!\!_{A+}(X)\simto W\!\!_{A+}(X)\simto\Om_{A+}(X),\\ \pp W\!\!_{A-}(X)\simto W\!\!_{A-}(X)\simto\Om_{A-}(X).\endaligned
\leqno(1.2.1)$$
\msn
{\it Proof.} The monoids in (\hl{1.2.1}{1.2.1}) are abelian groups by Proposition~\hl{P1.1b}{1.1b} (with $\F\eq\F'\eq\G\eq\G'$). We first show the isomorphism
\htt{1.2.2}{}
$$\pp W\!\!_{A+}(X)\simto\Om_{A+}(X).
\leqno(1.2.2)$$
Since the surjectivity follows from Proposition~\hl{P1.1a}{1.1a}, it is enough to show the injectivity.
Assume $(\F,S)$ is cobordant to $(\F',S')$ with $\F,\F'\in D^b_c(X,A)^{[0]}$. By Theorem~\hl{T1.1}{1.1}, there are self-dual $A$-complexes $(\F_i,S_i)$ for $i\in[0,2r]$ such that $(\F_0,S_0)\eq(\F,S)$, $(\F_{2r},S_{2r})\eq(\F',S')$, and $(\F_{i-1},S_{i-1})$ is directly subquotient cobordant to $(\F_i,S_i)$ with $\F_i\in D^b_c(X,A)^{[0]}$ for any $i\in[1,2r]$. Using Proposition~\hl{P1.1b}{1.1b}, we then see that the following self-dual complexes are directly cobordant to 0\,:
$$\aligned(\G_1,S_1):=\mopl_{i=1}^r&\bl((\F_i,S_i)\oplus(\F_i,-S_i)\br),\\(\G_2,S_2):=\mopl_{i=1}^r&\bl((\F_{i-1},S_{i-1})\oplus(\F_i,-S_i)\br).\endaligned$$
Since
$$(\F,S)\oplus(\G_1,S_1)=(\F',S')\oplus(\G_2,S_2),$$
we thus get the injectivity of (\hl{1.2.2}{1.2.2}).
\sk
It is now enough to show the surjectivity $\pp W\!\!_{A+}(X)\onto W\!\!_{A+}(X)$ for the proof of the first two isomorphisms of (\hl{1.2.1}{1.2.1}) (that is, for $+$); but this surjectivity follows from Propositions~\hl{P1.1a}{1.1a}--\hl{P1.1b}{b} considering
$$(\F,S)\oplus(\F,-S)\oplus\pp\Hc^0(\F,S).$$
The argument is similar with $+$ replaced by $-$.
This finishes the proof of Theorem~\hl{T1.2}{1.2}.
\par\htt{R1.2a}{}\msn
{\bf Remark~1.2a.} If $X$ is a point, $\pp W\!\!_{A+}({\rm pt})$ is identified with the classical Witt group $W(A)$, see \cite{MH} for the latter. (Recall that $A$ is a field of characteristic 0 in this subsection.)
It is also possible to verify that $\Om_{A+}({\rm pt})$ is isomorphic to $W(A)$ by an elementary argument as follows. (This is for a better understanding of the reader about cobordism relation.)
\sk
By Proposition~\hl{P1.1a}{1.1a}, any element of $\Om_{A+}({\rm pt})$ is represented by a non-degenerate symmetric bilinear form $S$ on a finite dimensional vector space $V$. By Theorem~\hl{T1.1}{1.1}, it is enough to show that, if $(V,S)$ is directly subquotient cobordant to $(V',S')$ with horizontal morphisms injective, then $(V',S')$ is the direct sum of $(V,S)$ and a metabolic form $(V'',S'')$ (which is a direct sum of hyperbolic forms).
\sk
We have the commutative diagram as in (\hl{1.1.4}{1.1.4})
$$\begin{array}{ccc}G_0V'&\into&V'\\ \rlap{\h{$\downarrow$}}\raise1mm\h{$\downarrow$}&\raise4mm\h{}\raise-2mm\h{}&\!\rlap{\h{$\downarrow$}}\raise1mm\h{$\downarrow$}\\ \Gr^G_0V'&\into&V'/G_{-1}V'\end{array}$$
associated to the direct subquotient cobordism relation. Here $G_i$ is the increasing filtration of $V'$ such that $\Gr^G_iV'=0$ ($i\ne -1,0,1$) and $\Gr^G_0V'=V$. Condition~(\hl{1.1.5}{1.1.5}) then implies the orthogonality
$$S'(G_0V',G_{-1}V')=0.$$
Taking an appropriate basis of $V'$, this means that $S'$ is represented by a matrix of the form
$$\begin{pmatrix}0&0&I\\0&S&B\\I&{}^t\!B&A\end{pmatrix}$$
Here $A$ is symmetric, $I$ is the identity matrix, and the matrix corresponding to $S$ is denoted also by $S$. Replacing the basis of $V'$ appropriately, we can then assume that $A\eq B\eq 0$ by a well-known argument in linear algebra (using ${}^t\!E_{\al,p,q}\1S'E_{\al,p,q}$, where the $(i,j)$-component of $E_{\al,p,q}\mi I$ is $\al$ if $(i,j)\eq(p,q)$, and it vanishes otherwise). So the desired isomorphism follows.
\par\htt{R1.2b}{}\msn
{\bf Remark~1.2b.} Using Theorems~\hl{T1.1}{1.1} and \hl{T1.2}{1.2}, we can construct the well-defined $L$-class transformation from $\Om_A(X)$ to $H_{2\1\ssb}(X,\Q)$ for {\it compact\1} algebraic varieties $X$ as in \cite{GM}, \cite{CS}, \cite{Ba1}, where a {\it cohomotopy} $f:X\to S^{2k}$ is used. Indeed, we see that the signature of the shifted restriction of $(\F,S)$ to the fiber of the north pole of $S^{2k}$ depends only on the cobordism class using the pushforward to the north pole. Although $f$ is not algebraic, we can apply a similar argument. Note also that $X$ may be replaced by a product with $S^{2m}$, and $(\F,S)$ by a shifted pull-back to the product as in \cite{Ba1}.
\par\htt{R1.2c}{}\msn
{\bf Remark~1.2c.} According to the referee, there are apparently only {\it surjections\1} from the cobordism groups and Witt groups in the literature (\cite{Bal}, \cite{BSY}, \cite{SW}, \cite{Wo}, \cite{Yo}, etc.) to those constructed in this paper, but they are actually {\it isomorphisms\1} using \cite[Thm.\,7.4]{Yo}, \cite[Prop.\,2.14]{SW}, \cite{Wo}, etc. This problem does not seem to be serious in case the reader is interested only in the $L$-classes. This issue may be better discussed in a different paper, since it seems rather technically complicated.
(Note that the cobordism relation is not used in the proof of Proposition~\hl{P1}{1}, since we show the independence of the {\it isomorphism class\1} of $(K_{\R}^{\ssb},S_{\R})$ under the choice of the polarization $S$ by constructing an {\it automorphism\1} of the complex $K_{\R}^{\ssb}$, see \hl{2.1}{2.1} below. This also applies to the corresponding result in \cite{FP}.)
\msn
{\bf 1.3.~Hodge modules.} For a smooth variety $X$ and for a subfield $A\subset\R$, we denote by ${\rm MF}_h(\D_X,A)$ the category of holonomic filtered $\D_X$-modules with $A$-structure $((M,F),K,\al)$. Here $(M,F)$ is a holonomic (analytic) $\D_X$-module with a good filtration $F$, $K$ is a $A$-complex with constructible cohomology sheaves on $X$, and the following isomorphisms of $\C$-complexes is given:
$$\al:{\rm DR}(M)\cong K\1{\otimes_A}\1\C,$$
see also Convention~(\hl{i}{i}) at the end of the introduction.
Recall that ${\rm DR}(M)$ denotes the shifted de Rham complex, see \cite{ypg}.
\sk
The category $\MH(X,A,w)$ of $A$-Hodge modules of weight $w$ on $X$ is a full subcategory of ${\rm MF}_h(\D_X,A)$. If $X$ is singular, this can be defined by using local embedding into smooth varieties. There is a strict support decomposition
\htt{1.3.1}{}
$$\MH(X,A,w)=\mopl_{Z\subset X}\,\MH_Z(X,A,w).
\leqno(1.3.1)$$
Here $Z$ runs over irreducible closed subvarieties of $X$, and an object of $\MH_Z(X,A,w)$ has {\it strict support\1} $Z$ (that is, its underlying $A$-complex is an intersection complex with local system coefficients).
Moreover we have the following.
\par\htt{T1.3}{}\msn
{\bf Theorem~1.3} (\cite[Thm.~2.21]{mhm}). {\it There is an equivalence of categories
\htt{1.3.2}{}
$$\MH_Z(X,A,w)={\rm VHS}^p_{\rm gen}(Z,A,w\mi d_Z),
\leqno(1.3.2)$$
where the right-hand side is the category of polarizable variations of $A$-Hodge structure of weight $w\mi d_Z$ defined on some Zariski-open subset of $Z$.
\sk
Moreover, the restriction of a polarization of Hodge module to the Zariski-open subset gives a polarization of variation of Hodge structure up to the sign $(-1)^{d_Z(d_Z-1)/2}$. Conversely, any polarization of variation of Hodge structure can be extended to a polarization of Hodge module up to the same sign.}
\sk
Note that variations of Hodge structure in this paper are assumed to have {\it unipotent\1} local monodromies.
\par\htt{R1.3}{}\msn
{\bf Remark~1.3.} A polarization of an $A$-Hodge module of weight $w$ is a $(-1)^w$-symmetric perfect pairing
\htt{1.3.3}{}
$$S:K{\otimes_A}K\to(\DD A_X)(-w).
\leqno(1.3.3)$$
There is a difference in signs of polarizations of Hodge modules and variations of Hodge structure as is explained above, see also \cite[5.2.12 and 5.4.1]{mhp} and Convention~(\hl{ii}{ii}) at the end of the introduction.
\par\htt{1.4}{}\msn
{\bf 1.4.~Mixed Hodge modules.} The category $\MHM(X,A)$ of mixed $A$-Hodge modules on a smooth variety $X$ is a full subcategory of ${\rm MF}_hW(\D_X,A)$, where the latter is the category of holonomic filtered $\D_X$-modules with $A$-structure $\M=\bl((M,F),K\br),\al)$, which are endowed with a finite filtration $W$ on $M,K$ in a compatible way with $\al$ so that
\htt{1.4.1}{}
$$\Gr^W_w\M\in\MH(X,A,w)\q(\forall\,w\in\Z).
\leqno(1.4.1)$$
However, the last condition is not enough, and it is rather complicated to give a precise definition for $\MHM(X,A)\subset{\rm MF}_hW(\D_X,A)$ (see \cite{def}) although we have the inclusion
\htt{1.4.2}{}
$$\MH(X,A,w)\subset\MHM(X,A)\q(\forall\,w\in\Z).
\leqno(1.4.2)$$
If $X$ is singular, it is defined by using local embedding into smooth varieties.
We can show that $\MHM(X,A)$ is an abelian category. Hence we have the derived category of bounded complexes of mixed Hodge modules $D^b\MHM(X,A)$ together with the standard cohomology functor
\htt{1.4.3}{}
$$H^{\ssb}:D^b\MHM(X,A)\to\MHM(X,A).
\leqno(1.4.3)$$
\sk
When $X\eq{\rm pt}$, there is an equivalence of categories
\htt{1.4.4}{}
$$\MHM({\rm pt},A)={\rm MHS}_A,
\leqno(1.4.4)$$
where the right-hand side is the category of graded-polarizable $A$-Hodge structures.
\sk
For a morphism of algebraic varieties $f:X\to Y$, we can construct the direct image functors
$$f_!,f_*:D^b\MHM(X,A)\to D^b\MHM(Y,A),$$
see \cite[Thm. 4.3]{mhm}, using Beilinson type resolution of $\M$ (see \cite{Be}).
We can show that the induced cohomological functor $H^jf_*:\MHM(X,A)\to\MHM(Y,A)$ is identified with the functor $\Hc^jf_*$ constructed in \cite{mhp} in the $f$ projective case.
\sk
The pull-back functors $f^*$, $f^!$ are defined respectively as left and right adjoint functors of the direct image functors $f_*$, $f_!$.
For Zariski-closed and open immersions $i:Y\into X$, $j:U\into X$ with $X\eq Y\,{\sqcup}\,U$, we have a distinguished triangle of functors
\htt{1.4.5}{}
$$j_!j^!\to{\rm id}\to i_*i^*\buildrel{[1]}\over\to.
\leqno(1.4.5)$$
This implies the well-definedness of the morphism $\Hdg$ in the introduction.
\sk
We have the following.
\par\htt{T1.4}{}\msn
{\bf Theorem~1.4} (\cite[Thm.~5.3.1]{mhp}). {\it Let $f:X\to Y$ be a projective morphism, and $\M\in\MH(X,A,w)$ with a polarization $S:K{\otimes_A}K\to(\DD A_X)(-w)$. Let $\ell\in H^2(X,A(1))$ be the first Chern class of a relative ample line bundle of $f$. Then}
\msn
\rlap{\rm(a)}\hskip.8cm\hangindent=.8cm
$H^jf_*\M\in\MH(Y,A,w{+}j)\q\q(j\in\Z),$
\msn
\rlap{\rm(b)}\hskip.8cm\hangindent=.8cm
$\ell^j:H^{-j}f_*\M\simto H^jf_*\M(j)\q\q(j>0),$
\msn
\rlap{\rm(c)}\hskip.8cm\hangindent=.8cm
{\it $(-1)^{j(j-1)/2}\1\RR f_*S\ssc(\ell^j{\otimes}{\rm id}):P_{\ell}\1\pp\Hc^{-j}f_*K{\otimes}P_{\ell}\1\pp\Hc^{-j}f_K\to(\DD A_Y)(j-w)$ is a polarization of $P_{\ell}\1H^{-j}f_*\M:={\rm Ker}\,\ell^{j+1}\subset H^{-j}f_*\M$ $\,\,(j\ges 0)\raise5mm\h{}$, where $P_{\ell}$ denotes the $\ell$-primitive part.}
\par\htt{R1.4}{}\msn
{\bf Remark~1.4.} Using the ``dictionary" in Convention~(\hl{ii}{ii}) at the end of the introduction, one can deduce the sign in the other sign system as in \cite{FP} from Theorem~\hl{T1.4}{1.4}\,(c), see also Remark~(iii) at the end of Section~1 in \cite{ypg} for a special case. In order to apply this ``dictionary", we may need also some equality like
$$(d{-}d')(d{-}d'{+}1)/2\pl(d{-}d')d'=d(d{-}1)/2\mi d'(d'{-}1)/2\pl(d{-}d'),$$
in the case of the {\it middle\1} primitive part. (We can restrict to this case, that is, $l\eq0$ in the notation of \cite{FP}, using the property of $\ep_m$ in (\hl{2.2.5}{2.2.5}) below.) Here $Z\subset Y$ is an irreducible closed subvariety of the target of a projective morphism $f:X\to Y$ with $X$ smooth, and $d\eq d_X$, $d'\eq d_Z$.
The left-hand side of the above equality comes from \cite{FP} (with $l\eq 0$), and the right-hand side from the dictionary (since $j(j{-}1)/2\eq0$ when $j\eq 0$).
Note that $d{-}d'$ coincides with the weight of the variation of Hodge structure associated to the direct factor of $H^0f_*(\Q_{h,X}[d_X])$ with strict support $Z$.
\bs\htt{S2}{}\bs
\vbox{\centerline{\bf 2. Proof of the main theorem and proposition}
\bsn
In this section we prove Theorems~\hl{T1}{1}--\hl{T3}{3} after showing Proposition~\hl{P1}{1}.}
\par\htt{2.1}{}\msn
{\bf 2.1.~Proof of Proposition~\hl{P1}{1}.} By Theorem~\hl{T1.3}{1.3}, it is enough to consider the case of a polarizable variation of Hodge structure of weight $w$ (using the functoriality of intermediate direct images).
Let $L$ be the underlying $\R$-local system of a polarizable variation of $\R$-Hodge structure $\V$ on a smooth complex algebraic variety $X$. Let
$$S:L\1{\otimes}\1L\to\R_X,\q S':L\1{\otimes}\1L\to\R_X$$
be polarizations of variation of Hodge structure. The polarizations $S,S'$ induce isomorphisms of variations of $\R$-Hodge structure
$$\V\buildrel{S}\over\cong(\DD\V)(-w)\buildrel{S'}\over\cong\V,$$
where $\DD\V$ denotes the dual variation of Hodge structure.
We thus get an endomorphism of variation of $\R$-Hodge structure, denoted by $\phi\in{\rm End}(\V)$, such that
\htt{2.1.1}{}
$$S'=S\ssc(\phi_{\R}\,{\otimes}\,{\rm id}),
\leqno(2.1.1)$$
where $\phi_{\R}\in{\rm End}(L)$ is the underlying morphism of $\phi$.
It is then enough to show that
\htt{2.1.2}{}
$$\h{$\phi_{\R}$ is semisimple and the eigenvalues of $\phi_{\R}$ are real numbers.}
\leqno(2.1.2)$$
Indeed, this implies the direct sum decomposition
\htt{2.1.3}{}
$$L=\mopl_i\,L_i\q\h{with}\q L_i={\rm Ker}(\phi_{\R}\mi\al_i)\subset L,
\leqno(2.1.3)$$
with $\al_i$ the eigenvalues of $\phi_{\R}$ (in particular, $\phi_{\R}|_{L_i}=\al_i$). We see that the $\al_i$ must be {\it positive\1} considering linear combinations of the two polarizations with positive coefficients. (Note that polarizations are stable by sum and also by multiplication by positive numbers.) The two polarizations are then identified by using the positive square roots of the eigenvalues.
\sk
For the proof of (\hl{2.1.2}{2.1.2}), we can restrict to a point of $X$ so that we have a polarizable real Hodge structure $H=(H_{\R},(H_{\C},F))$ with two polarizations $S$, $S'$ and an endomorphism of real Hodge structure $\phi$ such that
\htt{2.1.4}{}
$$S'(u,v)=S(\phi u,v)\q(u,v\in H_{\R}).
\leqno(2.1.4)$$
Since the Weil operator $C$ can be defined on $H_{\R}$, and $S,S'$ are $(-1)^w$-symmetric with $w$ the weight of Hodge structure, we have a {\it positive-definite\1} symmetric form
$$S_C:H_{\R}\times H_{\R}\to\R$$
defined by
$$S_C(u,v)=S(u,Cv),$$
(see \cite{De1}) and similarly for $S'_C$. Here the symmetry of $S_C$ follows from
\htt{2.1.5}{}
$$S(u,Cv)=(-1)^wS(Cu,v)=S(v,Cu).
\leqno(2.1.5)$$
It is then enough to show that $\phi$ is symmetric with respect to $S_C$, that is,
\htt{2.1.6}{}
$$S_C(\phi u,v)=S_C(u,\phi v),
\leqno(2.1.6)$$
(since $\phi$ can be expressed by a symmetric matrix taking an orthonormal basis of $H_{\R}$). But this symmetry follows from
\htt{2.1.7}{}
$$S(\phi u,Cv)=S'(u,Cv)=S'(v,Cu)=S(\phi v,Cu)=S(u,C\phi v).
\leqno(2.1.7)$$
So (\hl{2.1.2}{2.1.2}) is proved. This finishes the proof of Proposition~\hl{P1}{1}.
\par\htt{R2.1}{}\msn
{\bf Remark~2.1.} The following assertion is shown in a recent version of \cite{FP}: Let $L$ be a $\R$-local system on a complex (or even topological) manifold $U$. Let $\{S_t\}_{t\in I}$ be a smooth 1-parameter family of non-degenerate symmetric (or skew-symmetric) pairing of $L$, where $I$ is an open interval. Then there is a smooth 1-parameter family of automorphisms $A_t\in{\rm Aut}(L)$ locally on $I$ such that the pairings $S_t\ssc(A_t\1{\otimes}\1 A_t)$ are independent of $t$. 
Here one does {\it not\1} assume that $L$ is {\it semisimple.} So the situation is quite different from the case of Proposition~\hl{P1}{1}, and the $A_t$ are {\it not necessarily semisimple.} The proof uses integral curves of time dependent vector fields.
In the case where $U\eq\C^*$ and the monodromy of $L$ is unipotent with only one Jordan block, one can verify the assertion directly, where the pairing $S_{t_0}$ for some fixed point $t_0\in I$ is represented by an {\it anti-diagonal\1} matrix, choosing a base point of $U$.
\msn
{\bf 2.2.~Proof of Theorem~\hl{T3}{3}.} By (\hl{1.4.1}{1.4.1}--\hl{1.4.2}{2}), it is enough to consider the case of pure Hodge modules with strict support. (Note that any short exact sequence of mixed Hodge modules induces short exact sequences by taking the $\Gr^W_i$, and pure Hodge modules are semisimple by Theorem~\hl{T1.3}{1.3}.) The well-definedness of the morphism $\Pol$ follows from Proposition~\hl{P1}{1}.
\sk
We first show the compatibility of $\Pol$ with the pushforward by {\it projective\1} morphisms. Let $f:X\to Y$ be a projective morphism of complex varieties. Let $\M\in\MH_Z(X,A,w)$ with $K$ the underlying $\R$-complex and $S$ a polarization. Put
$$\ep_m:=(-1)^{m(m+1)/2}\q\q(m\in\Z).$$
We have
\htt{2.2.1}{}
$$\aligned\Pol(f_*[\M])&=\msum_j\,(-1)^j\1\Pol[H^jf_*\M],\\ f_*(\Pol[\M])&=\ep_w[(\RR f_*K,\RR f_*S)]=\ep_w[(\pp\Hc^0\RR f_*K,\pp\Hc^0\RR f_*S)],\endaligned
\leqno(2.2.1)$$
where the last equality follows from Proposition~\hl{P1.1a}{1.1a}.
\sk
By the hard Lefschetz property (see Theorem~\hl{T1.4}{1.4}\,(b)), we have the Lefschetz decomposition
\htt{2.2.2}{}
$$f_*[\M]=\msum_{j\ges0}\,\msum_{k=0}^j\,(-1)^j[P_{\ell}H^{-j}f_*\M(-k)],
\leqno(2.2.2)$$
where $P_{\ell}H^{-j}f_*\M:={\rm Ker}\,\ell^{j+1}\subset H^{-j}f_*\M$ is the primitive part for $j\ges 0$ (and it is 0 otherwise).
Note that $(-k)$ for $k\in\Z$ is the Tate twist which changes the weight by $2k$, and $\ell$ is a morphism of degree 2.
\sk
We denote the primitive part of the underlying $\R$-complex by
$$P_{\ell}^{-j}:=P_{\ell}\pp\Hc^{-j}\1\RR f_*K\q\q(j\ges 0).$$
This has the self-pairing defined by $\RR f_*S\ssc(\ell^j{\otimes}{\rm id})$, which is denoted by $f_*S$ for simplicity. It has also a self-pairing defined by a polarization of Hodge module, which is denoted by $S_p$. By Theorem~\hl{T1.4}{1.4}\,(c), we have
\htt{2.2.3}{}
$$[(P_{\ell}^{-j},f_*S)]=(-1)^{j(j-1)/2}[(P_{\ell}^{-j},S_p)].
\leqno(2.2.3)$$
By the Lefschetz decomposition (\hl{2.2.2}{2.2.2}) together with Theorem~\hl{T1.4}{1.4}\,(a), we get that
\htt{2.2.4}{}
$$\Pol(f_*[\M])=\msum_{j\ges0}\,\msum_{k=0}^j\,(-1)^j\1\ep_{w-j+2k}\1[(P^{-j}_{\ell},S_p)].
\leqno(2.2.4)$$
\sk
Observe that the restriction of $\ep_m$ to even or odd numbers has the {\it alternating\1} property:
\htt{2.2.5}{}
$$\ep_m=(-1)^i\,\,\,\h{if}\,\,\,\,m\eq 2i\,\,\,\,\h{or}\,\,\,\,m\eq 2i{-}1.
\leqno(2.2.5)$$
This implies that the summation over $k\in[0,j]$ in (\hl{2.2.4}{2.2.4}) is an alternating sum. Hence it vanishes if $j$ is odd, and the summation over $k\in[1,j]$ vanishes if $j$ is even. Using (\hl{2.2.3}{2.2.3}), we then get that
\htt{2.2.6}{}
$$\aligned\Pol(f_*[\M])&=\ep_w\,\msum_{k\ges0}\,(-1)^k\1[(P^{-2k}_{\ell},S_p)]\\&=\ep_w\,\msum_{k\ges0}\,[(P^{-2k}_{\ell},f_*S)],\endaligned
\leqno(2.2.6)$$
since $\ep_{w-2k}=\ep_w(-1)^k$ and $(-1)^{i(i-1)/2}=(-1)^k$ if $i\eq 2k$. Thus the compatibility with the pushforward by projective morphisms is proved. (Note that $\ell^k$ is a morphism of even degree, and does not produce a sign.)
\sk
We now consider the case $f$ is non-projective. We may assume that $\M$ is a pure Hodge module with strict support $X$ (in particular, $X$ is irreducible). We proceed by induction on $d_X\,(=\dim X)$. The assertion holds if $d_X=0$. Assume $d_X\eq d$ and the commutativity holds when $d_X<d$.
There is a projective birational morphism $\pi:X'\to X$ such that $g:=f\ssc\pi:X'\to Y$ is also projective. Consider the diagram
$$\begin{array}{ccccc}K_0(\MHM(X',A))&\buildrel{\pi_*}\over\longrightarrow&K_0(\MHM(X,A))&\buildrel{f_*}\over\longrightarrow&K_0(\MHM(Y,A))\\ \,\,\,\downarrow\!{\scriptstyle\rm Pol}&\raise5mm\h{}\raise-2mm\h{}&\,\,\,\downarrow\!{\scriptstyle\rm Pol}&&\,\,\,\downarrow\!{\scriptstyle\rm Pol}\\ \Om_{\R}(X')&\buildrel{\pi_*}\over\longrightarrow&\Om_{\R}(X)&\buildrel{f_*}\over\longrightarrow&\Om_{\R}(Y)
\end{array}$$
(This is partially commutative as is explained in (\hl{2.2.8}{2.2.8}) below.)
\sk
Since $\pi$ is birational, there are $\M'\in\MH_{X'}(X',A,w)$ and $\M''\in\MH(Z,A,w)$ such that
\htt{2.2.7}{}
$$H^0\pi_*\M'=\M\oplus(i_Z)_*\M''\q\h{in}\,\,\,\MH(X,A,w),
\leqno(2.2.7)$$
by Theorems~\hl{T1.3}{1.3} and \hl{T1.4}{1.4}\,(a), where $Z\buildrel{i_Z}\over\into X$ is a proper closed subvariety.
\sk
Since $g\eq f\ssc\pi$ and $\pi$ are projective, we have the following equalities in $\Om_{\R}(Y)$\,:
\htt{2.2.8}{}
$$\Pol(f_*(\pi_*[\M']))=f_*(\pi_*(\Pol[\M']))=f_*(\Pol(\pi_*[\M'])).
\leqno(2.2.8)$$
By (\hl{2.2.7}{2.2.7}), there is $\xi\in K_0(\MHM(Z,A))$ such that
\htt{2.2.9}{}
$$\pi_*[\M']=[\M]\pl(i_Z)_*\1\xi\q\h{in}\,\,\,K_0(\MHM(X,A)),
\leqno(2.2.9)$$
replacing $Z\subset X$ so that ${\rm Supp}\,H^j\pi_*\M'\subset Z$ for $j\ne 0$. (Here the decomposition theorem is not needed, since we consider the equality in the Grothendieck group.) By inductive hypothesis, we have
\htt{2.2.10}{}
$$\Pol(f_*(i_Z)_*\1\xi)=f_*(\Pol((i_Z)_*\1\xi))\q\h{in}\,\,\,\Om_{\R}(Y).
\leqno(2.2.10)$$
Then (\hl{2.2.8}{2.2.8}--\hl{2.2.10}{10}) imply the quality
\htt{2.2.11}{}
$$\Pol(f_*([\M]))=f_*(\Pol([\M]))\q\h{in}\,\,\,\Om_{\R}(Y).
\leqno(2.2.11)$$
This completes the proof of Theorem~\hl{T3}{3}.
\par\htt{2.3}{}\msn
{\bf 2.3.~Proofs of Theorems~\hl{T1}{1} and \hl{T2}{2}.} Let $C$ be the mapping cone of the morphism~(\hl{8}{8}) in the introduction. Its image in the Grothendieck group of mixed Hodge modules $K_0({\MHM}(X,\R))$ comes from $\eta_x\in K_0({\MHM}(\{x\},\R))$ for $x\in\Si$. The Hodge signature of $\eta_x$ is given by $\sit_x$, and vanishes by assumption. Hence its image in the cobordism group of $X$ vanishes using the compatibility with the pushforward by the inclusion $\Si\into X$. (Recall that an element of the classical Witt group over $\R$ is determined only by the signature, see \cite{MH} and also Remark~\hl{R1.2a}{1.2a} for the relation with the cobordism group.) Theorem~\hl{T1}{1} thus follows.
\sk
The argument is similar for Theorem~\hl{T2}{2}. The converse holds in this case, since the difference is contained in the image of $H_0(\Si,\Q)$ in $H_{2\ssb}(X,\Q)$, and is given by $\sum_{x\in\Si}\sit_x$ (using the connectivity of $X$). This finishes the proofs of Theorems~\hl{T1}{1} and \hl{T2}{2}.
\par\htt{R2.3a}{}\msn
{\bf Remark~2.3a.} The assumption of Theorem~\hl{T1}{1} is satisfied if the $d_{X_i}$ are {\it even\1} and $X$ has only {\it isolated hypersurface\1} singularities defined analytic-locally by holomorphic functions $f$ with {\it semisimple\1} Milnor monodromies, for instance, if $X$ is locally isomorphic to a cone of a smooth projective hypersurface of odd dimension, that is, $f$ is a homogeneous polynomial, or more generally, $f$ is a semi-weighted-homogeneous polynomial $f\eq\sum_{\alpha\gess 1}\,f_{\alpha}$ where the $f_{\alpha}$ are weighted homogeneous polynomials of degree $\alpha\gess 1$ with weights of variables fixed, and $f_1$ has an isolated singularity at 0. This is a $\mu$-constant deformation of $f_1$, and has semisimple Milnor monodromy.
\sk
Indeed, assuming $X$ has an isolated hypersurface singularity defined by $f$, set
$$K_{f,x}:={\rm Ker}\bl(N:H^{d_X}(F_{\!f},\Q)_1\to H^{d_X}(F_{\!f},\Q)_1(-1)\br).$$
Here $H^{d_X}(F_{\!f},\Q)_1$ denotes the unipotent monodromy part of the vanishing cohomology with $F_{\!f}$ the Milnor fiber, and $N:=\log T_u$ with $T=T_sT_u$ the Jordan decomposition of the monodromy (and $d_X:=\dim X$). It is quite well-known that
\htt{2.3.1}{}
$$H_{X,x}^j:=\Hc^j({\rm IC}_X\Q[-d_X])_x=\begin{cases}K_{f,x}(1)&(j=d_X{-}1),\\ \,0&(j\ne 0,\,d_X{-}1).\end{cases}
\leqno(2.3.1)$$
Note that $H^{d_X}(F_{\!f},\Q)_1$ is {\it pure\1} with weight $d_X{+}1$, if the Milnor monodromy is {\it semisimple\1} (that is, $N\eq 0$), see \cite{St} and also \cite[(5.1.6.2)]{mhp}.
\sk
It is rather easy to show the $W$-graded quotients of (\hl{2.3.1}{2.3.1}), which is sufficient for our purpose, by using the short exact sequence of mixed Hodge modules
\htt{2.3.2}{}
$$0\to\Q_{h,X}[d_X]\to\psi_{f,1}\Q_{h,Y}[d_X]\buildrel{\rm can}\over\longrightarrow\varphi_{f,1}\Q_{h,Y}[d_X]\to 0.
\leqno(2.3.2)$$
Here $f$ is defined on an ambient complex manifold $Y$ containing locally $X$, and $\psi_{f,1}$, $\varphi_{f,1}$ denote the unipotent monodromy part of $\psi_f$, $\varphi_f$ respectively.
The weight filtration $W$ on $\psi_{f,1}\Q_{h,Y}$, $\varphi_{f,1}\Q_{h,Y}$ is the monodromy filtration shifted by $d_X$ and $d_X{+}1$ respectively. So there are isomorphisms
$$\aligned N^j:\Gr^W_{d_X+j}\psi_{f,1}(\Q_Y[d_X])&\simto\Gr^W_{d_X-j}\psi_{f,1}(\Q_Y[d_X])(-j),\\ N^j:\Gr^W_{d_X+1+j}\varphi_{f,1}(\Q_Y[d_X])&\simto\Gr^W_{d_X+1-j}\varphi_{f,1}(\Q_Y[d_X])(-j),\endaligned$$
for $j>0$ together with the $N$-{\it primitive decomposition\1}
$$\Gr^W_k(\psi_{f,1}\Q_Y[d_X])=\mopl_{j\ges 0}\,N^j\h{$_P$}\Gr^W_{k+2j}(\psi_{f,1}\Q_Y[d_X])(j),$$
for $k\in\Z$, and similarly with $\psi_{f,1}$ replaced by $\varphi_{f,1}$. Here the $N$-{\it primitive part\1} is defined by
$$\aligned{}_P\Gr^W_{d_X+j}(\psi_{f,1}\Q_Y[d_X])&:={\rm Ker}\,N^{j+1}\subset\Gr^W_{d_X+j}(\psi_{f,1}\Q_Y[d_X]),\\{}_P\Gr^W_{d_X+1+j}(\varphi_{f,1}\Q_Y[d_X])&:={\rm Ker}\,N^{j+1}\subset\Gr^W_{d_X+1+j}(\varphi_{f,1}\Q_Y[d_X]),\endaligned$$
for $j\ges 0$, and these vanish otherwise. By construction $N$ coincides with the composition of 
$$\psi_{f,1}\Q_Y[d_X]\buildrel{\rm can}\over\longrightarrow\varphi_{f,1}\Q_Y[d_X]\buildrel{\rm Var}\over\longrightarrow\psi_{f,1}\Q_Y(-1)[d_X].$$
Here $\,{\rm can}\,$ is {\it surjective,} $\,{\rm Var}\,$ is {\it injective,} hence $\varphi_{f,1}\Q_Y[d_X]$ is identified with ${\rm Coim}\,N$, and $_h$ is sometimes omitted to simplify the notation. Note that $\M\mapsto\Gr^W_k\M$ is an exact functor of mixed Hodge modules, and commutes with ${\rm Ker}$, ${\rm Coim}$.
\sk
For $j>0$, we can then deduce the isomorphisms
$$\aligned{}_P\Gr^W_{d_X+j}(\psi_{f,1}\Q_Y[d_X])&={}_P\Gr^W_{d_X+j}(\varphi_{f,1}\Q_Y[d_X]),\\ \Gr^W_{d_X-j}(\Q_X[d_X])&={\rm Ker}\,N\subset\Gr^W_{d_X-j}\psi_{f,1}(\Q_Y[d_X])\\&={\rm Ker}\,N\subset\Gr^W_{d_X-j+2}\varphi_{f,1}(\Q_Y[d_X])(1),\endaligned$$
using the $N$-primitive decompositions for $\psi_{f,1}$, $\varphi_{f,1}$ together with (\hl{2.3.2}{2.3.2}).
\sk
For $j=0$, it is known (see for instance \cite[(4.5.9)]{mhm}) that
\htt{2.3.3}{}
$$\Gr^W_{d_X}(\Q_{h,X}[d_X])={\rm IC}_{h,X}\Q.
\leqno(2.3.3)$$
So the $W$-graded quotients of (\hl{2.3.1}{2.3.1}) follow by using the short exact sequence
$$0\to W_{d_X-1}(\Q_{h,X}[d_X])\to\Q_{h,X}[d_X]\to{\rm IC}_{h,X}\Q\to 0.$$
Note that $H^{d_X}(F_{\!f},\Q)_1=\varphi_{f,1}\Q_{h,Y}[d_X]$.
We have the vanishing for $0<j<d_X{-}1$ in (\hl{2.3.1}{2.3.1}), since $\Q_{h,X}[d_X]$ is a mixed Hodge module.
(The isomorphism (\hl{2.3.1}{2.3.1}) without $\Gr^W$ can be shown for instance using the octahedral diagram for the composition ${\rm Var}\ssc{\rm can}=N$.)
\par\htt{R2.3b}{}\msn
{\bf Remark~2.3b.} A typical example with {\it semisimple\1} Milnor monodromy is an ordinary double point defined locally by $f\eq\msum_{i=0}^{d_X}\,x_i^2\,+$ higher terms. In the $d_X$ {\it odd\1} case, the reduced modified Euler-Hodge signature is given by $\sit_x\eq(-1)^{(d_X-1)/2}$, since
$$\h{$H_{X,0}^j\eq\Q\bl(-\tfrac{j}{2}\br)$ if $j\eq 0,d_X{-}1$, and 0 otherwise.}$$
\sk
As an example with {\it non-semisimple\1} Milnor monodromy, we have
$$f=\msum_{i=1}^4\,x_i^{p_i}+x_1x_2x_3x_4+\h{higher terms}\q\q\bl(\1\msum_{i=1}^4\,\tfrac{1}{p_i}<1\br),$$
with $d_X\eq3$. Here ``higher terms" are with respect to the {\it Newton polyhedron.} If the $p_i$ are {\it mutually prime,} the unipotent monodromy part of the Milnor fiber cohomology has dimension 3, and we have $H_{X,x}^2\eq\Q$ (where the unipotent monodromy part of the Milnor fiber cohomology has only one Jordan block with size 3), hence $\sit_x\eq1$. This can be verified by using \cite{JKSY}, \cite{des}. The latter shows that there is a Jordan block of size 3 in this case. The spectrum is given by setting $\gamma_{\langle p\rangle}(t):=\msum_{i=1}^{p-1}\,t^{i/p}$ as follows.
\htt{2.3.4}{}
$${\rm Sp}_f(t)=\msum_{k=1}^3\,t^k+\msum_{i=1}^4\,\gamma_{\langle p_i\rangle}(t)(t{+}t^2)+\msum_{i<j}\,\gamma_{\langle p_i\rangle}(t)\gamma_{\langle p_j\rangle}(t)\,t.
\leqno(2.3.4)$$
One can also employ a computer program like Singular \cite{DGPS}, for instance, in the case $\{p_i\}\eq\{3,4,5,7\}$ as below (where the spectral numbers are shifted by $-1$).
\ms
\vbox{\small\sf\verb#LIB "gmssing.lib";#
\sk
\verb#ring R = 0, (x,y,z,w), ds; poly f=x^3+y^4+z^5+w^7+x*y*z*w;#
\sk
\verb#spprint(spectrum(f));#
\sk
\verb#spprint(sppairs(f));#}
\msn
(Here {\small\sf spectrum(f)} gives multiplicities, but {\small\sf sppairs(f)} does not.) These can be used to construct an example of a globally irreducible variety with $\msum_x\,\sit_x\eq0$ using Remark~\hl{R2.3c}{2.3c} below.
\par\htt{R2.3c}{}\msn
{\bf Remark~2.3c.} Let $f,g\in\C[x_1,\dots,x_n]$ be Newton non-degenerate polynomials with isolated singularities at the origin. Let $F(z_0,\dots,z_n),G(z_0,\dots,z_n)$ be homogeneous polynomials of the same degree $d$ such that
$$F(1,x_1,\dots,x_n)\eq f(x_1,\dots,x_n),\q G(x_1,\dots,x_n,1)\eq g(x_1,\dots,x_n).$$
Put
\vskip-7mm
$$H:=F+c_nG+\msum_{i=1}^{n-1}\,c_iz_i^d.$$
Then the projective hypersurface $X$ defined by $H$ has isolated singularities at the two points defined by $z_i\eq0$ ($i\,{\ne}\,0$) and ($i\,{\ne}\,n$) respectively, and they have the same non-degenerate Newton boundaries as $f,g$ unless the degree $d$ is very small. Moreover, if the $c_i$ are sufficiently general, $X$ is non-singular outside the two points. (The last assertion may be skipped by using resolutions of singularities.)
\par\htt{R2.3d}{}\msn
{\bf Remark~2.3d.} It is also possible construct an example of a globally irreducible variety having a singular point $x$ with $\sit_x\eq 0$ by a cancellation. For this, consider two 4-dimensional subvarieties $Z_1,Z_2\subset\PP^8$ defined respectively by
$$\aligned x_1^2\pl 3x_2^2\pl\msum_{i=1}^3\,y_i^2\eq 4x_0^2,&\q z_i\eq0\,\,(i\in[1,3]),\\ 3x_1^2\pl x_2^2\pl\msum_{i=1}^3\,z_i^2\eq 4x_0^2,&\q y_i\eq0\,\,(i\in[1,3]),\endaligned$$
with $x_0,x_1,x_2,y_1,y_2,y_3,z_1,z_2,z_3$ projective coordinates of $\PP^8$. These subvarieties intersect at the following 4 points:
$$\bl\{x_1\eq{\pm}\1x_0,\,x_2\eq{\pm}\1x_0,\,y_i\eq z_i\eq0\,\,(i\in[1,3])\br\}.$$
There is an involution $\iota$ of $\PP^8$ such that $\iota(Z_1)\eq Z_2$. This is defined by
$$\iota:(x_0,x_1,x_2,y_1,y_2,y_3,z_1,z_2,z_3)\mapsto(x_0,x_2,x_1,z_1,z_2,z_3,y_1,y_2,y_3).$$
(This describes the corresponding involution of $\C^9$, and similarly for the expressions below.) We have the morphism $\pi_1:Z_1\to\PP^4$ defined by
$$\pi_1:(x_0,x_1,x_2,y_1,y_2,y_3,z_1,z_2,z_3)\mapsto(x_0,x_2,y_1,y_2,y_3).$$
Set $\pi_2:=\pi_1\ssc\iota:Z_2\to\PP^4$. These are quotient morphisms under the actions defined by multiplication by $-1$ on $x_1$ and $x_2$ respectively. We may assume that the images in $\PP^4$ of the above 4 intersection points are the two singular points of the hypersurface $X\sst\PP^4$ constructed in Remarks~\hl{R2.3b}{2.3b}-\hl{R2.3c}{c} (using an automorphism of $\PP^4$). The desired irreducible variety is the quotient of $\pi_1^{-1}(X)\cup\pi_2^{-1}(X)$ by the action of $\iota$. (Note that the $\pi_i^{-1}(X)$ are irreducible.) We can eliminate unnecessary singularities at the fixed points of $\iota$ using resolutions of singularities. So the quotient of the singularities at the non-fixed points only remains.
\par\htt{R2.3e}{}\msn
{\bf Remark~2.3e.} If $X$ is the affine cone of a {\it smooth projective\1} variety $Z$ with $0$ the vertex of cone, we have the isomorphism
\htt{2.3.5}{}
$$H_{X,0}^j:=\Hc^j({\rm IC}_X\Q[-d_X])_0=H^j_{\rm prim}(Z,\Q)\q(j\in\Z),
\leqno(2.3.5)$$
where the right-hand side denotes the primitive cohomology of $Z$. This follows from the Thom-Gysin sequence, see for instance \cite{RSW}. Here $Z$ is not assume to be a projective {\it hypersurface.} This can be used to construct some other examples. (The construction does not seem necessarily trivial. For instance, one can take a surjection from the projective cone to projective space of the same dimension, and take the base change by a double cover of the projective space ramified along a divisor which does not contain the image of the vertex of the projective cone and has an appropriate isolated singularity.)
\msn
{\bf 2.4.~Explicit calculation.} Let $X$ be a compact threefold having only one singular point $P$, which is an ordinary $m$-ple point with $m\ges 2$ (that is, $X$ is locally defined by a semi-homogeneous polynomial $f=\sum_{k\ges m}f_k$ where the $f_k$ are homogeneous polynomials of degree $k$ and $f_m$ has an isolated singularity at 0). Let $\pi:\Xt\to X$ be the blow-up along $P$. Then its exceptional divisor is a smooth projective surface $Z\subset\PP^3$ of degree $m$ which is defined by $f_m$. There is the equality in $K_0({\rm Var}/X)$\,:
\htt{2.4.1}{}
$$[X]=[\Xt]+[P]-[Z].
\leqno(2.4.1)$$
We have
\vskip-7mm
$$\pp\Hc^j(\RR\pi_*\Q_{\Xt}[3])\cong\begin{cases}{\rm IC}_X\Q&(j\eq0),\\ \Q_P&(j\eq{\pm}1),\\ \,0&(\h{otherwise}),\end{cases}$$
where Tate twists are omitted.
By Proposition~\hl{P1.1a}{1.1a}, this implies the equality
$$[(\RR\pi_*\Q_{\Xt}[3],S)]=[({\rm IC}_X\Q,S)]\q\h{in}\,\,\,\,\Om(X).$$
The image of $[X]$ in $\Om(X)$ is then expressed as
\htt{2.4.2}{}
$$[({\rm IC}_X\Q,S)]+[(\Q_P,S)]-[(H^2(Z,\Q),S)].
\leqno(2.4.2)$$
Here the contribution of the surface $Z$ to $\Om(X)$ is given by $H^2(Z)$ using Proposition~\hl{P1.1a}{1.1a}. Its non-primitive part cancels with $[(\Q_P,S)]$ after the scalar extension by $\Q\into\R$. However, the primitive part $H^2_{\rm prim}(Z,\Q)$ remains. This gives a difference between the homology $L$-class which is the specialization of the Hirzebruch class and the intersection complex $L$-class of $X$ (as is seen from the proof of Theorem~\hl{T2}{2}).
\par\htt{R2.4}{}\msn
{\bf Remark~2.4.} If we consider the case where $d_X$ is {\it even\1} and $X$ has only singular points which are ordinary $m$-ple points as above (where $m$ may be depend on the points), then we can see that the two classes coincide as in Theorem~\hl{T2}{2}. Indeed, the difference is given by the sum of $[H^{d_Z}(Z,\R),S]$ with $Z$ exceptional divisors (here we can use also Remarks~\hl{R2.3a}{2.3a}--\hl{R2.3b}{b} together with a $\mu$-constant deformation), and $S$ is {\it skew-symmetric,} since $d_Z\eq d_X\mi 1$ is {\it odd.}
\par\htt{2.5}{}\msn
{\bf 2.5.~Proof of Proposition~\hl{P2}{2}.} Let $X$ be a compact surface having only one singular point $P$, which is an ordinary double point (that is, $A_1$-singularity). Let $\pi:\Xt\to X$ be the blow-up along $P$ with $C\,(=\PP^1)\subset\Xt$ the exceptional divisor. We have the equality in $K_0({\rm Var}/X)\,$:
\vskip-4mm\htt{2.5.1}{}
$$[X]=[\Xt]+[P]-[\PP^1].
\leqno(2.5.1)$$
Here $\PP^1$ does not contribute to the class of $X$ in $\Om(X)$ (since $H^1(\PP^1)=0$). We have the canonical isomorphism
\vskip-6mm
$$\RR\pi_*\Q_{\Xt}[2]={\rm IC}_X\Q\oplus\Q_P,$$
which gives the decomposition in $\Om(X)$\,:
$$[(\RR\pi_*\Q_{\Xt}[2],S)]=[({\rm IC}_X\Q,S)]+[(\Q_P,S')],$$
where $S$ on the left-hand side and on the intersection complex is the canonical pairing, and $S'$ in the last term is induced from $S$ on the left-hand side. The image of $[X]$ in $\Om(X)$ is then given by
\htt{2.5.2}{}
$$[({\rm IC}_X\Q,S)]+[(\Q_P,S')]+[(\Q_P,S)],
\leqno(2.5.2)$$
where the last $S$ is the canonical pairing. However, $[(\Q_P,S')]$ does not cancel with $[(\Q_P,S)]$ in $\Om(X)$, since $S'$ is described as
\htt{2.5.3}{}
$$\Q_P\otimes\Q_P\ni(u,v)\mapsto-2\1uv\in\Q.
\leqno(2.5.3)$$
Note that the self-intersection number of $C\subset\Xt$ is $-2$ as is well-known.
\sk
This implies a counter-example to Theorem~\hl{T1.3}{1.3} with $\Q$-coefficients (using the pushforward by $X\to{\rm pt}$). Indeed, the symmetric cobordism group $\Om_{\Q+}({\rm pt})$ is identified with the classical Witt group $W(\Q)$ by Theorem~\hl{T1.2}{1.2}. Moreover the image of $[\langle-2\rangle_{\Q_2}]$ by the morphism
$$\psi^1:W(\Q_2)\to W(\FF_2)$$
is equal to $[\langle -1\rangle_{\FF_2}]\eq[\langle 1\rangle_{\FF_2}]$, which does not vanish, see Remarks~\hl{R2.5a}{2.5a}--\hl{R2.5b}{b} below. (One may use also $\psi^0:W(\Q_3)\to W(\FF_3)$, since $\langle -2\rangle_{\FF_3}\eq\langle 1\rangle_{\FF_3}$ and $2\1[\langle 1\rangle_{\FF_3}]\ne 0$.)
Here $\langle\al\rangle_K$ with $K$ a field denotes the bilinear form on the 1-dimensional vector space $K$ defined by $\al\in K$. This finishes the proof of Proposition~\hl{P2}{2}.
\par\htt{R2.5a}{}\msn
{\bf Remark~2.5a.} For a prime number $p$, there are two natural morphisms of abelian groups
\htt{2.5.4}{}
$$\psi^k:W(\Q_p)\to W(\FF_p)\q\q(k\eq 0,1),
\leqno(2.5.4)$$
sending $[\langle up^i\rangle_{\Q_p}]$ to $[\langle\overline{u}\rangle_{\FF_p}]$ if $i\equiv k$ mod 2, and to 0 otherwise. Here $u\in\Z_p$ is invertible, and $\overline{u}$ is its image in $\FF_p^{\1*}$, see \cite[IV.1.2]{MH}.
\par\htt{R2.5b}{}\msn
{\bf Remark~2.5b.} There are isomorphisms of abelian groups (see \cite[IV.1.5]{MH})
\htt{2.5.5}{}
$$W(\FF_p)\cong\begin{cases}\Z/2\1\Z&(p\eq 2),\\ \Z/4\1\Z&(p\equiv 3\,\,\,{\rm mod}\,\,4),\\ \Z/2\1\Z\times\Z/2\1\Z&(p\equiv 1\,\,\,{\rm mod}\,\,4).\end{cases}
\leqno(2.5.5)$$


\begin{thebibliography}{CMSS2}
\bibitem[Bal]{Bal} Balmer, P., Witt groups, in Handbook of $K$-theory, Springer, Berlin, 2005, pp.~539--576.
\bibitem[Ban1]{Ba1} Banagl, M., Topological Invariants of Stratified Spaces, Springer, Berlin, 2007.
\bibitem[Ban2]{Ba2} Banagl, M., Topological and Hodge $L$-classes of singular covering spaces and varieties with trivial canonical class, Geometriae Dedicata 199 (2019), 189--224.
\bibitem[BFM]{BFM} Baum, P., Fulton, W., MacPherson, R.D., Riemann-Roch for singular varieties, Publ.\ Math.\ IHES 45 (1975), 101--145.
\bibitem[Be]{Be} Beilinson, A., On the derived category of perverse sheaves, Lect.\ Notes in Math.\ 1289, Springer, Berlin (1987), 27--41.
\bibitem[BBD]{BBD} Beilinson, A., Bernstein, J., Deligne, P., Faisceaux pervers, Ast\'erisque 100, Soc.\ Math.\ France, Paris, 1982.
\bibitem[Bi]{Bi} Bittner, F., The universal Euler characteristic for varieties of characteristic zero, Compos.\ Math.\ 140 (2004), 1011--1032.
\bibitem[BM]{BM} Borho, W., MacPherson, R., Partial resolutions of nilpotent varieties, Ast\'erisque 101-102 (1983), 23--74.
\bibitem[BSY]{BSY} Brasselet, J.-P., Sch\"urmann, J., Yokura, S., Hirzebruch classes and motivic Chern classes for singular spaces, J. Topol.\ Anal.\ 2 (2010), 1--55.
\bibitem[BS]{BS} Brasselet, J.-P., Schwartz, M.-H., Sur les classes de Chern d'une ensemble analytique complexe, in Caract\'eristice d'Euler-Poincar\'e, S\'eminaire E.N.S. 1978-1979. Ast\'erisque 82-83 (1981), 93--148.
\bibitem[Br]{Br} Brieskorn, E., Die Monodromie der isolierten Singularit\"aten von Hyperfl\"achen, Manuscripta Math. 2 (1970), 103--161.
\bibitem[CH]{CH} Calm\`es, B., Hornbostel, J., Tensor-triangulated categories and dualities, Theory Appl.\ Categ.\ 22 (2009), 136–198 and 200.
\bibitem[CS]{CS} Cappell, S.E., Shaneson,J.L., Stratifiable maps and topological invariants, J.\ Amer.\ Math.\ Soc.\ 4 (1991), 521--551.
\bibitem[CMSS1]{CMSS1} Cappell, S.E., Maxim, L.G., Sch\"urmann, J., Shaneson, J.L., Characteristic classes of complex hypersurfaces, Adv.\ Math.\ 225 (2010), 2616-- 2647.
\bibitem[CMSS2]{CMSS2} Cappell, S.E., Maxim, L.G., Sch\"urmann, J., Shaneson, J.L., Equivariant characteristic classes of singular complex algebraic varieties, Comm.\ Pure Appl.\ Math.\ 65 (2012), 1722--1769.
\bibitem[DGPS]{DGPS} Decker, W., Greuel, G.-M., Pfister, G., Sch\"onemann, H., {\sc Singular} 4.2.0 --- A computer algebra system for polynomial computations, available at http://www.singular.uni-kl.de (2020).
\bibitem[De1]{De1} Deligne, P., The\'orie de Hodge II, Publ.\ Math. IHES 40 (1971), 5--58.
\bibitem[De2]{De2} Deligne, P., Cohomologie \`a supports propres, SGA4 XVII, Springer Lect.\ Notes in Math.\ 305 (1973), 250-461.
\bibitem[De3]{De3} Deligne, P., Positivit\'e: signe, I (manuscrit, 16-2-84), II (manuscrit, 6-11-85).
\bibitem[FP]{FP} Fern\'andez de Bobadilla, J., Pallar\'es, I., The Brasselet-Sch\"urmann-Yokura conjecture on $L$-classes of projective varieties (arxiv:2007.11537).
\bibitem[GM]{GM} Goresky, M., Macpherson, R.D., Intersection homology theory, Topology 19 (1980), 135--162.
\bibitem[Hi]{Hi} Hirzebruch, F., Topological Method in Algebraic Geometry, Springer, Berlin, 1966.
\bibitem[JKSY]{JKSY} Jung, S.-J., Kim, I.-K., Saito, M., Yoon, Y., Spectrum of non-degenerate functions with simplicial Newton polyhedra (arxiv:1911.09465).
\bibitem[Ma]{Ma} MacPherson, R.D., Chern classes for singular algebraic varieties, Ann.\ Math.\ 100 (1974), 423--432.
\bibitem[MS]{MS} Maxim, L.G., Sch\"urmann, J., Characteristic classes of singular toric varieties, Comm.\ Pure Appl.\ Math.\ 68 (2015), 2177--2236.
\bibitem[MSS]{MSS} Maxim, L.G., Saito, M., Sch\"urmann, J., Symmetric products of mixed Hodge modules, J.\ Math.\ Pures Appl.\ (9) 96 (2011), 462--483.
\bibitem[MH]{MH} Milnor, J., Husem\"oller, D., Symmetric Bilinear Forms, Springer, Berlin, 1973.
\bibitem[RSW]{RSW} Reichelt, T., Saito, M., Walther, U., Dependence of Lyubeznik numbers of cones of projective schemes on projective embeddings, Selecta Math.\ (N.S.) 27 (2021), Paper No.\ 6, 22 pp.
\bibitem[Sa1]{mhp} Saito, M., Modules de Hodge polarisables, Publ.\ RIMS, Kyoto Univ.\ 24 (1988), 849--995.
\bibitem[Sa2]{mhm} Saito, M., Mixed Hodge modules, Publ.\ RIMS, Kyoto Univ.\ 26 (1990), 221--333.
\bibitem[Sa3]{mhc} Saito, M., Mixed Hodge complexes on algebraic varieties, Math.\ Ann.\ 316 (2000), 283--331.
\bibitem[Sa4]{def} Saito, M., On the definition of mixed Hodge modules (arxiv:1307.2140).
\bibitem[Sa5]{ypg} Saito, M., A young person's guide to mixed Hodge modules, in Hodge theory and $L^2$-analysis, Adv.\ Lect.\ Math.\ 39, Int.\ Press, Somerville, MA, 2017, 517--553.
\bibitem[Sa6]{lw} Saito, M., Weight zero part of the first cohomology of complex algebraic varieties (arXiv:1804.03632).
\bibitem[Sa7]{des} Saito, M., Descent of nearby cycle formula for Newton non-degenerate functions (arxiv:2004.12367).
\bibitem[Sch]{Schu} Sch\"urmann, J., Characteristic classes of mixed Hodge modules, in Topology of Stratified Spaces, MSRI Publ.\ 58, 2011, pp.~419--470.
\bibitem[SW]{SW} Sch\"urmann, J., Woolf, J., Witt groups of abelian categories and perverse sheaves, Ann.\ $K$-theory 4 (2019), 621--670.
\bibitem[Schw]{Schw} Schwartz, M.-H., Classes caract\'eristiques d\'efinies par une stratification d'une vari\'et\'e analytique complexe, C.\ R.\ Acad.\ Sci.\ Paris 260 (1965), 3262--3264 and 3535--3537.
\bibitem[St]{St} Steenbrink, J.H.M., Mixed Hodge structure on the vanishing cohomology, in Real and complex singularities, Sijthoff and Noordhoff, Alphen aan den Rijn, 1977, pp. 525--563.
\bibitem[Ve]{Ve} Verdier, J.-L., Cat\'egories d\'eriv\'ees, Etat 0. in SGA 4 1/2, Springer Lect.\ Notes in Math.\ 569 (1977), 262--308.
\bibitem[Wo]{Wo} Woolf, J., Witt groups of sheaves on topological spaces, Comment.\ Math.\ Helv.\ 83 (2008), 289--326.
\bibitem[Yo]{Yo} Youssin, B., Witt groups of derived categories, $K$-theory 11 (1997), 373--395.
\end{thebibliography}
\end{document}